\documentclass[12pt,twoside]{amsart}
\usepackage[latin1]{inputenc}
\usepackage{amsmath,amsfonts,amssymb,amscd,amsthm,amsxtra}

\usepackage{a4wide}
\usepackage{url}
\usepackage{enumerate,verbatim}
\usepackage{ifthen}

\usepackage{url,textcomp}
\newcommand{\no}{\textnumero}

\usepackage{xr1}

\usepackage{pstricks}
\newcommand{\ddate}[1]{}

\newcommand{\noteQ}[2]{#2}
\newcommand{\ppaper}[4][ ]{}
\newcommand{\paper}[4][ ]{}

\numberwithin{equation}{section}

\makeatletter
  \DeclareRobustCommand\em
         {\@nomath\em \ifdim \fontdimen\@ne\font >\z@
                       \upshape \else \slshape \fi}
\makeatother

\newtheoremstyle{mythm}
  {9pt}
  {9pt}
  {\slshape}
  {0pt}
  {\bfseries}
  {.}
  { }
  {\thmname{#1} \thmnumber{#2}\thmnote{ (#3)}}

\theoremstyle{mythm} 
\newtheorem{Theorem}{Theorem}[section]
\newtheorem{Lemma}[Theorem]{Lemma}

\newtheorem{Proposition}[Theorem]{Proposition}

\newtheorem{Corollary}[Theorem]{Corollary}

\theoremstyle{definition} 
\newtheorem{Definition}[Theorem]{Definition}

\newtheorem{Remark}[Theorem]{Remark}

\newtheorem{Example}[Theorem]{Example}
 
\newcommand{\file}[1]{{\url{#1}}}

\newcommand{\IC}{\mathbf{C}}                     
\newcommand{\IM}{I\!\! M}                        
\newcommand{\IN}{\mathbf{N}}                     
\newcommand{\IZ}{\mathbf{Z}}                     
\newcommand{\SG}{\mathfrak{S}}                   
\renewcommand{\phi}{{\varphi}}                   
\newcommand{\ox}{\otimes}                        
\newcommand{\NCPartitions}{NC}                     
\newcommand{\NC}{\NCPartitions}                    
\newcommand{\abs}[1]{\left\lvert #1 \right\rvert}  

\newcommand{\IncAlg}{\mathfrak{I}}                 
\newcommand{\exchm}{\alg{E}} 
\newcommand{\exchF}{{\alg{F}}} 
\newcommand{\exchCF}{{\mathrm{CF}}} 
\newcommand{\exchGrad}{{\mathrm{GT}}} 
\newcommand{\exchFA}{{\alg{F}_a}} 
\newcommand{\exchB}{{\mathrm{Boo}}} 
\newcommand{\exch}{\mbox{$\exchm$}} 
\DeclareMathOperator*{\arrowprod}{\overset{\rightarrow}{\prod}}    

\newlength{\tmpl}

\newcommand{\bub}[1]{{\overset{\circ}{#1}}}      

\newcommand{\indep}{\perp\kern-3pt\perp}

\newcommand{\alg}[1]{{\mathcal{#1}}}                

\DeclareMathOperator{\IE}{\mathbf{E}}                     
\DeclareMathOperator{\tr}{{\mathrm tr}}                     

\begin{document}

\title[Cumulants in Noncommutative Probability Theory I]%
{Cumulants in Noncommutative Probability Theory I.
  Noncommutative Exchangeability Systems}
\author{Franz Lehner}

\thanks{Supported by the European Network \no{}HPRN-CT-2000-00116
and the Austrian Science Fund (FWF), Project \no{}R2-MAT}

\address{
Franz Lehner\\
In\-sti\-tut f\"ur Mathe\-ma\-tik C\\
Tech\-ni\-sche Uni\-ver\-si\-t\"at Graz\\
Stey\-rer\-gas\-se 30, A-8010 Graz\\
Austria}
\email{lehner@finanz.math.tu-graz.ac.at}


\date{\today}

\begin{abstract}
  Cumulants linearize convolution of measures.
  We use a formula of Good to define noncommutative cumulants in a very general setting.
  It turns out that the essential property needed is exchangeability
  of random variables. 
  Roughly speaking the formula says that cumulants are moments of 
  a certain ``discrete Fourier transform'' of a random variable.
  This provides a simple unified method to understand the
  known examples of cumulants, like classical, free and various $q$-cumulants.
\end{abstract}

\maketitle{}

\tableofcontents{}
\goodbreak{}
\hfill\parbox{.5\textwidth}{{\itshape
It will be shown that the formulae are much simplified by the use of
cumulative moment functions, or semi-invariants, in place of the crude moments.}
\par
R.A.~Fisher \cite{Fisher:1929:moments}}

\vskip2em

The object of this series of papers is a unified treatment of cumulants.
A wide variety of cumulants has been defined in different contexts,
like classical cumulants and free cumulants, the latter being the most well-known
noncommutative example. Each of these examples is tailored for a
certain notion of independence, but all of them share a certain similarity.
It will turn out that this is no coincidence and that all these definitions
have a common source, namely a certain exchangeability relation.
This rather general condition will be the starting point for our definition
of independence.

There have been axiomatic approaches to noncommutative independence,
for example in the work of Sch\"urmann (see, e.g., \cite{Schurmann:1995:direct,BenGhorbalSchurmann:2001:algebraic})
in the context of co- and bialgebras. The axioms there, while natural, are quite
rigorous and it was shown by Speicher \cite{Speicher:1997:universal}
that under these axioms there are only three possibilities -- 
classical, free and boolean independence.

In another vein, there were attempts to adapt classical cumulants
to noncommutative situations, cf.\ Hegerfeldt  \cite{Hegerfeldt:1985:noncommutative}.
These considerations are however confined to tensor product constructions.

The aim of the present paper is to show that certain combinatorial aspects of independence
hold in the context of exchangeability.
It may be disputed if the term ``independence'' is justified here.
There are certain combinatorial analogies
with the notion of independence of classical random variables,
notably visible in part~II (\cite{Lehner:2002:Cumulants2}),
while other properties fail.
The main drawback in this setting is that the joint distribution of independent noncommutative
random variables is not determined by the distributions of the individual random variables.
This is one of the main axioms in Sch\"urmann's approach and already seen
to fail for $q$-independence, see \cite{vanLeeuwenMaassen:1996:obstruction}.
As a consequence our notion of independence is non-constructive,
that is, an infinite family of interchangeable algebras must be given \emph{a priori}.
(An exception to this is fermionic independence
(section~\ref{sec:GeneralCumulantExamples:Fermions} below)
where the presence of additional structure,
namely a $\IZ_2$-grading, provides for another invariant and independent algebras can
be constructed by means of graded tensor products.)
If one accepts these drawbacks there still remains a rich unified combinatorial theory
comprising many known examples and opening the field for new ones.

The paper roughly splits into two halves.

In
sections~\ref{sec:GeneralCumulants:Introduction}--\ref{sec:GeneralCumulants:BasicTransformations}
we use a formula of Good to define cumulants and  ``independence''
 with respect to so-called \emph{exchangeability systems}.
The basic properties of cumulants are almost immediately obvious from this formula.
Alternatively, after expanding Good's formula and collecting equal terms
one rediscovers the well-known definition of cumulants via M\"obius inversion on the lattice
of set partitions in full generality.
From a computational point of view, the second definition is more efficient and a
large number of combinatorial formulas from classical statistics can be
transferred to the general setting.

In section~\ref{sec:GeneralCumulants:Examples} 
we use the general machinery to recompute several
known examples of cumulants and exhibit why a particular kind of cumulants is
the ``right'' one for a certain notion of independence. 

In subsequent papers \cite{Lehner:2002:Cumulants2,Lehner:2002:Cumulants3}
we will treat characterizations of so-called generalized Gaussian random
variables (or \emph{generalized Brownian motions}) 
and exchangeable random variables arising from Fock space constructions.

\emph{Acknowledgements.}
The author is grateful to Dan~Voiculescu for a discussion, during which he suggested the
notion of ``exchangeability system'', which makes the concepts much clearer.
We also acknowledge the comments of two anonymous referees on an earlier version
of this paper.


\section{Introduction and definitions}
\label{sec:GeneralCumulants:Introduction}

\subsection{Classical Cumulants}
Cumulants were introduced by Thiele in his 1889 book 
under the name of \emph{semi-invariants}, but entered
the wider scene of statistics only with Fisher's fundamental paper  \cite{Fisher:1929:moments}
under the name of \emph{cumulative moment functions}. Shortly afterwards,
the name \emph{cumulants} was commonly adopted.
We refer to \cite{Mattner:1999:what} for the analytical aspects of classical cumulants
and to \cite{Hald:2000:early} for their history.
Here the focus will be on the combinatorial aspects of cumulants.
\begin{Definition}
  \label{def:GeneralCumulants:ClassicalCumulants}
  Let $X$ be a random variable with moments~$m_n=m_n(X)$ and denote
  $$
  \mathcal{F}_X(z) = \IE e^{zX} = \sum_{n=0}^\infty \frac{m_n}{n!}\,z^n
  $$
  its formal \emph{Fourier-Laplace transform} or
  \emph{exponential moment generating function}, considered as a formal power series.
  The coefficients $\kappa_n=\kappa_n(X)$ of its formal logarithm
  $$
  \log \mathcal{F}_X(z) = \sum_{n=1}^\infty \frac{\kappa_n}{n!}\,z^n
  $$
  are called the \emph{(classical) cumulants} of~$X$.
\end{Definition}

Equivalently, classical cumulants can be defined by the recursion formula

\begin{equation}
  \label{eq:GeneralCumulants:ClassicalCumulantRecursion}
  \kappa_n = m_n - \sum_{k=1}^{n-1} \binom{n-1}{k-1} \kappa_k m_{n-k}
  .
\end{equation}

In this paper we consider cumulants for \emph{noncommutative} or \emph{quantum probability}
spaces.
\begin{Definition}
  A noncommutative probability space is a pair $(\alg{A},\phi)$ of
  a complex unital algebra $\alg{A}$ equipped with a unital linear functional $\phi$,
  which is called the \emph{expectation}. The elements of $\alg{A}$ are called
  \emph{(noncommutative) random variables}.
  Usually $\alg{A}$ will be a $C^*$-algebra and $\phi$ a faithful state.
  More generally, an \emph{operator-valued} noncommutative probability space
  is a unital algebra $\alg{A}$ together with a unital subalgebra $\alg{B}$
  and a \emph{conditional expectation} $\psi:\alg{A}\to\alg{B}$,
  i.e.,
  a linear map $\psi$ which satisfies the identity
  $\psi(bab')=b\,\psi(a)\,b'$ for all $a\in\alg{A}$ and $b$,~$b'\in\alg{B}$.
  Such an algebra is also called \emph{$\alg{B}$-valued probability space}
  and its elements are \emph{$\alg{B}$-valued random variables}.
\end{Definition}
In order to define cumulants, one needs a notion of \emph{independence}
or, as it will turn out, \emph{exchangeability}.
The most prominent example of independence in noncommutative probability
is Voiculescu's free probability theory \cite{VDN:1992:free}.
Many concepts from classical probability have analogues in free probability,
among them are cumulants.
Existence of free cumulants was already proved in \cite{Voiculescu:1985:symmetries},
and a beautiful systematic theory was developed by R.~Speicher 
\cite{Speicher:1994:multiplicative} with many applications.

Another notion of cumulants (``partial cumulants'') was introduced even earlier
by von Waldenfels \cite{vonWaldenfels:1973:approach,vonWaldenfels:1975:interval}
and turned out to be connected to boolean independence 
\cite{SpeicherWoroudi:1997:boolean,BozejkoSpeicher:1991:psi}
associated to Bozejko's ``regular'' free product 
of states \cite{Bozejko:1987:uniformly}.
Other kinds of cumulants appear throughout noncommutative probability theory
and will be reviewed in section~\ref{sec:GeneralCumulants:Examples}.

The common characteristics of these cumulants can be summarized in the following 
properties, which in the classical case can easily be deduced from
Definition~\ref{def:GeneralCumulants:ClassicalCumulants}.
To any random variable $X$ having moments $m_n(X)$ of all orders,
there is associated a sequence $K_n(X)$ with the following properties.
\begin{subequations}
  \label{eq:GeneralCumulants:CumulantProperties1-3}
  \begin{enumerate}
   \item Additivity. If $X$ and $Y$ are independent random variables, then
    \begin{equation}
      \label{eq:GeneralCumulants:axiom1}
      K_n(X+Y) = K_n(X)+K_n(Y).
    \end{equation}
   \item Homogeneity. For any scalar~$\lambda$ the $n$-th cumulant is $n$-homogeneous:
    \begin{equation}
      \label{eq:GeneralCumulants:axiom2}
      K_n(\lambda X) = \lambda^n K_n(X).
    \end{equation}
   \item There exists a polynomial $P_n$ in $n-1$ variables 
    without constant term such that
    \begin{equation}
      \label{eq:GeneralCumulants:axiom3}
      m_n(X) = K_n(X) + P_n(K_1(X),K_2(X),\dots,K_{n-1}(X))
      .
    \end{equation}
  \end{enumerate}
\end{subequations}
  ``Independence'' here means classical (resp.\ free, boolean) independence
  in the case of classical (resp.\ free, boolean) cumulants.

\subsection{Good's formula}

The aim of this paper is to define cumulants in a uniform way.
In section~\ref{ssec:independence} we introduce an appropriate notion of independence
which is based on exchangeability.
The axioms are satisfied by all known examples, which are reviewed in
section~\ref{sec:GeneralCumulants:Examples}.
In the future we hope to give new examples.
Our definition is based on a formula of Good \cite{Good:1975:new} for classical cumulants,
which shows up as a curiosity in the exercise sections of some textbooks of
statistics. While it is less useful in classical statistics, where much more
powerful methods of Fourier analysis are available, 
it will turn out to be very useful in noncommutative situations.

\begin{Theorem}[{Good \cite{Good:1975:new}}]
  Let $X$ be a random variable and $X^{(k)}$, $k=1,2,\dots,n$ be i.i.d.\ copies of $X$.
  Let $\omega$ be a primitive $n$-th root of unity and set
  $$
  X^{\omega} = \omega X^{(1)} + \omega^2 X^{(2)} + \dots + \omega^n X^{(n)}
  $$
  Then
  \begin{equation}
    \label{eq:GeneralCumulants:GoodFormula1}
    \kappa_n(X) = \frac{1}{n}\,\IE [(X^\omega)^n]
  \end{equation}
\end{Theorem}
The original proof consisted of two pages of computations,
but it was realized shortly afterwards that there is a three line proof
\cite{Good:1977:new}, 
based on the properties 
\eqref{eq:GeneralCumulants:CumulantProperties1-3}
together with a simple symmetry consideration.
For the reader's convenience we include this proof here.

\begin{proof}
  We evaluate the right hand side of \eqref{eq:GeneralCumulants:GoodFormula1}
  using property \eqref{eq:GeneralCumulants:axiom3} above:
  $$
  m_n(X^\omega)
  = \kappa_n(X^\omega) 
    + P_n(\kappa_1(X^\omega),\kappa_2(X^\omega),\dots,\kappa_{n-1}(X^\omega))
  ;
  $$
  now the cumulants of $X^\omega$ can be evaluated using properties
  \eqref{eq:GeneralCumulants:axiom1} and \eqref{eq:GeneralCumulants:axiom2}:
  \begin{align*}
    \kappa_m(\sum_k \omega^k X^{(k)})
    &= \sum_k \omega^{km} \kappa_m(X^{(k)})\\
    &= \sum_k \omega^{km} \kappa_m(X)
  \end{align*}
  and
  $$
  \sum_k \omega^{km} = 
  \begin{cases}
    n & \text{if $n$ divides $m$}\\
    0 & \text{otherwise}
  \end{cases}
  $$
  In particular the cumulant vanishes for $m<n$ and since $P_n$ has no constant term,
  the only contribution comes from $\kappa_n(X^\omega)$.
\end{proof}

\subsection{Posets and M\"obius inversion}
\label{sec:GeneralCumulant1:MoebiusInversion}
There is an alternative approach to cumulants using M\"obius inversion on the lattice
of set partitions.
The M\"obius function of a poset was introduced in a systematic manner by Rota
\cite{Rota:1964:foundationsI,DoubiletRotaStanley:1972:foundationsVI}.
Let $(P,\leq)$ be a (finite) partially ordered set, in short a \emph{poset}.
The incidence algebra $\IncAlg(P)=\IncAlg(P,\IC)$ is the algebra of 
functions supported on the set of pairs $\{(x,y)\in P\times P : x,y\in P; x\leq y\}$
with convolution
$$
f*g(x,y) = \sum_{x\leq z\leq y} f(x,z)\, g(z,y)
$$
For example, if $P$ is the $n$-set~$\{1,2,\dots,n\}$ with the natural order,
then $\IncAlg(P)$
is the algebra of $n\times n$ upper triangular matrices.
In general the algebra~$\IncAlg(P)$ has the identity $\delta(x,y)$ and a function $f\in\IncAlg(P)$
is invertible if and only if $f(x,x)$ is invertible for every $x\in P$.
The function
$\zeta(x,y)\equiv 1$ is called \emph{Zeta function}.
It is invertible and
its inverse is called the \emph{M\"obius function} of~$P$, 
denoted $\mu(x,y)$.
For functions $F,G:P\to\IC$ we have the fundamental equivalence
 (``M\"obius inversion formula'')
$$
\left(
  \forall x\in P:  F(x) = \sum_{y\leq x} G(y)
\right)
\qquad\iff\qquad
\left(
  \forall x\in P:  G(x) = \sum_{y\leq x} F(y)\,\mu(y,x)
\right)
$$
The poset $P$ is a \emph{lattice} if supremum and infimum operations exist.

\subsection{Partitions}

We will be working with the lattice of set partitions $\Pi_n$ and some of its sublattices.
\begin{Definition}
  A \emph{partition} of a set $S$ is a set $\pi=\{\pi_1,\pi_2,\dots,\pi_k\}$
  of pairwise disjoint nonempty subsets of~$S$ such that $\bigcup \pi_j=S$.
  Equivalently, a partition of $S$ corresponds to an equivalence relation $\sim_\pi$ on $S$
  where $i\sim_\pi j$ if $i$ and $j$ lie in the same block.
  The components $\pi_j$ of $\pi$ will be referred to as \emph{blocks} or \emph{classes}
  of $\pi$.
  The set of partitions of a set $S$ will be denoted by $\Pi_S$,
  or, if $S=[n]:=\{1,2,\dots,n\}$, we will abbreviate it as $\Pi_n$.
  It forms a lattice under the \emph{refinement order},
  where $\pi\leq\sigma$ if every block of $\pi$ is contained in some block of $\sigma$.
  In this ordering there is a maximal element $\hat1_n$ consisting of only one
  block and a minimal element $\hat0_n$ consisting of $n$ singletons.
\end{Definition}

Partitions can be visualized by diagrams, where the points are drawn on a line
and those points which lie in one block are connected by an arc.

For the examples in section~\ref{sec:GeneralCumulants:Examples}
we will be interested in various classes of partitions of the $n$-set
$[n]$ in which the order on $[n]$ will be important.
\begin{Definition}
  \label{def:GeneralCumulants:partitions}
  \begin{enumerate}
   \item A partition $\pi\in\Pi_n$ is \emph{noncrossing} if there is no quadruple
    of elements $i<j<k<l$ s.t.\ $i\sim_\pi k$, $j\sim_\pi l$ and $i\not\sim_\pi j$.
    The noncrossing partitions of order $n$ form a lattice which we denote by $\NC_n$.
   \item A block~$B$ of a noncrossing partition $\pi$ is \emph{inner} 
    if there are elements $i$, $j\not\in B$ 
    such that $i<k<j$ for all $k\in B$ and $i\sim_\pi j$.
    The other blocks are called the \emph{outer blocks} of $\pi$.
   \item An \emph{interval partition} is a partition $\pi$ for which every block is an interval.
    Equivalently, this means that~$\pi$ is noncrossing and
    all blocks of $\pi$ are outer.
   \item A partition $\pi$ is \emph{connected} if the picture of $\pi$ is a connected graph.
    The \emph{connected components} of $\pi$ are the maximal connected subpartitions of $\pi$.
   \item The \emph{noncrossing closure} of a partition $\pi$ is the smallest noncrossing
    partition which dominates $\pi$.
   \item A partition $\pi\in\Pi_n$ is \emph{irreducible} if the elements~$1$ and~$n$
    are in the same connected component.
    Every partition $\pi$ can be ``factored'' into irreducible factors.
   \item The \emph{interval closure} of a partition $\pi$ is the smallest interval
    partition which dominates $\pi$.
  \end{enumerate}
  Different types of partitions
  are shown in figure~\ref{fig:GeneralCumulants:Partitions}.
%
%
\begin{center}
  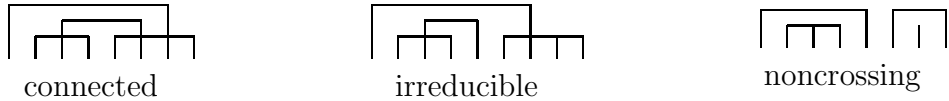
\begin{figure}[htbp]
    \begin{minipage}{.3\textwidth}
      \begin{center}
        \begin{picture}(80,20.4)(1,0)
          \put(10,0){\line(0,1){20.4}}
          \put(20,0){\line(0,1){8.4}}
          \put(30,0){\line(0,1){14.4}}
          \put(40,0){\line(0,1){8.4}}
          \put(50,0){\line(0,1){8.4}}
          \put(60,0){\line(0,1){14.4}}
          \put(70,0){\line(0,1){20.4}}
          \put(80,0){\line(0,1){8.4}}
          \put(20,8.4){\line(1,0){20}}
          \put(30,14.4){\line(1,0){30}}
          \put(10,20.4){\line(1,0){60}}
          \put(50,8.4){\line(1,0){30}}
        \end{picture}
        
        connected
      \end{center}
    \end{minipage}
%
%
    \begin{minipage}{.3\textwidth}
      \begin{center}
      \begin{picture}(90,20.4)(1,0)
        \put(10,0){\line(0,1){20.4}}
        \put(20,0){\line(0,1){8.4}}
        \put(30,0){\line(0,1){14.4}}
        \put(40,0){\line(0,1){8.4}}
        \put(50,0){\line(0,1){14.4}}
        \put(60,0){\line(0,1){8.4}}
        \put(70,0){\line(0,1){20.4}}
        \put(80,0){\line(0,1){8.4}}
        \put(90,0){\line(0,1){8.4}}
        \put(20,8.4){\line(1,0){20}}
        \put(30,14.4){\line(1,0){20}}
        \put(10,20.4){\line(1,0){60}}
        \put(60,8.4){\line(1,0){30}}
      \end{picture}

        irreducible
      \end{center}
    \end{minipage}
%
%
    \begin{minipage}{.3\textwidth}
      \begin{center}
        \begin{picture}(80,14.4)(1,0)
          \put(10,0){\line(0,1){14.4}}
          \put(20,0){\line(0,1){8.4}}
          \put(30,0){\line(0,1){8.4}}
          \put(40,0){\line(0,1){8.4}}
          \put(50,0){\line(0,1){14.4}}
          \put(60,0){\line(0,1){14.4}}
          \put(70,0){\line(0,1){8.4}}
          \put(80,0){\line(0,1){14.4}}
          \put(20,8.4){\line(1,0){20}}
          \put(10,14.4){\line(1,0){40}}
          \put(70,8.4){\line(1,0){0}}
          \put(60,14.4){\line(1,0){20}}
        \end{picture}
        
        noncrossing
      \end{center}
    \end{minipage}
    \label{fig:GeneralCumulants:Partitions}
    \caption{Typical partitions}
  \end{figure}
\end{center}
\end{Definition}

\begin{Proposition}
  \label{prop:GeneralCumulants:NCclosure}
  The noncrossing closure of a partition $\pi\in\Pi_n$ is obtained by
  putting all elements of each connected component into one block.
  Consequently a partition is noncrossing if and only if 
  every connected component consists of exactly one block.
  Another consequence is the fact that a partition $\pi$
  is irreducible if and only if $1\sim_{\hat\pi} n$ in its noncrossing closure
  $\hat\pi$.

  Similarly, the interval closure of $\pi$ is obtained by putting
  the elements of each irreducible factor into one block.
  Consequently an interval partition is characterized by the
  property that each irreducible factor consists of exactly one block.
\end{Proposition}

Many formulas in this paper will involve
partitions induced by index sequences
and for these partitions the following notation will be convenient.

\begin{Definition}
  \label{def:GeneralCumulants:kernel}
  Let~$f:[n]\to X$ be a function from the $n$-set $[n]$ to
  some set  $X$.
  The partition $\pi\in\Pi_n$ corresponding to the equivalence
  relation $i\sim_\pi j$ $\iff$ $f(i)=f(j)$ is called the
  \emph{kernel} of~$f$ and denoted~$\ker f$.
\end{Definition}

\subsection{Exchangeability and independence}
\label{ssec:independence}

There is a variant of Good's formula for multivariate cumulants
which will serve as a definition in the following situation.
\begin{Definition}
  \label{def:GeneralCumulants:independence}
  Let $(\alg{A},\phi)$ be a noncommutative probability space.
  An \emph{exchangeability system} $\alg{E}$ for $(\alg{A},\phi)$
  consists of a noncommutative probability space $(\alg{U},\tilde\phi)$
  and an infinite family $\alg{J}=(\iota_k)_{k\in\IN}$
  of state-preserving embeddings 
  $\iota_k:\alg{A}\to \alg{A}_k\subseteq \alg{U}$, 
  which we conveniently denote by $X\mapsto X^{(k)}$,
  such that the algebras $\alg{A}_j$ are \emph{interchangeable}
  with respect to $\tilde\phi$:
  for any family~$X_1,X_2,\dots,X_n\in \alg{A}$,
  and for any choice of indices~$i_1,i_2,\dots,i_n$
  the expectation is invariant under any permutation~$\sigma\in\SG_\infty$
  in the sense that
  $$
  \tilde\phi(X_1^{(i_1)} X_2^{(i_2)} \cdots X_n^{(i_n)})
  =
  \tilde\phi(X_1^{(\sigma(i_1))} X_2^{(\sigma(i_2))} \cdots X_n^{(\sigma(i_n))})
  .
  $$
  In other words, the value of the expectation only depends on
  the kernel of the map $h:j\mapsto i_j$, i.e.,
  the partition $\pi$ of $\{1,2,\dots,n\}$ made up from the equivalence classes
  of the equivalence   relation $j\sim_\pi k$ $\iff$ $i_j=i_k$.
  We will denote this value by $\phi^{\alg{E}}_\pi(X_1,X_2,\dots,X_n)$
  or $\phi_\pi(X_1,X_2,\dots,X_n)$ if the choice of $\alg{E}$ is
  clear from context. 

  Similarly, for a subset~$B\subseteq \{1,2,\dots,n\}$ (resp., a partition~$\pi$ of
  a subset~$B$)
  we will abbreviate the expectation
  $
  \phi_B(T_1,T_2,\dots,T_n) = \phi(\prod_{j\in B} T_j)
  $
  (ordered product) and $\phi_\pi(T_1,T_2,\dots,T_n) = \phi_\pi(T_j:j\in B)$.

  We will say that subalgebras $\alg{B},\alg{C}\subseteq \alg{A}$ are
  \emph{$\alg{E}$-exchangeable} or, more suggestively, \emph{$\alg{E}$-independent}
  if for any choice of random variables~$X_1,X_2,\dots,X_n \in \alg{B}\cup\alg{C}$
  and subsets~$I,J\subseteq \{1,\dots,n\}$ such that~$I\cap J=\emptyset$,
  $I\cup J=\{1,\dots,n\}$, $X_i\in \alg{B}$ for $i\in I$  and $X_i\in \alg{C}$ for $i\in J$,
  we have the identity
  $$
  \phi_\pi(X_1,X_2,\dots,X_n)
  = \phi_{\pi'}(X_1,X_2,\dots,X_n)
  $$
  whenever~$\pi$, $\pi'\in\Pi_n$ are partitions with $\pi|_I=\pi'|_I$ and $\pi|_J=\pi'|_J$.
  We say that two families of random variables~$(X_i)_{i\in I}$ and~$(Y_j)_{j\in J}$ are
  $\alg{E}$-exchangeable if the algebras they generate have this
  property.
\end{Definition}
\begin{Remark}
  \label{rem:GeneralCumulants:Independence}
  In other words, \exch-independence means that
  if~$\rho=\{I_\alg{B},I_\alg{C}\}$ is a partition as above
  then for any map~$h:[n]\to\IN$ the expectation
  $$
  \phi(X_1^{(h(1))} X_2^{(h(2))} \dotsm X_n^{(h(n))})
  $$
  is unchanged if we modify $h$ in such a way that the partition $\rho\wedge\ker h$
  does not change.
  
  Also note that for a given sequence $X_1,\dots,X_n$ there
  may be different choices for~$I_\alg{B}$ and~$I_\alg{C}$, if
  some of the $X_i$ lie in the intersection $\alg{B}\cup\alg{C}$.
\end{Remark}

\begin{Example}
  As an example, assume that the subalgebras $\alg{B}$ and $\alg{C}\subseteq\alg{A}$
  are \exch{}-independent in the above sense, then any noncommutative polynomial
  $P(X_1,\dots,X_n,Y_1,\dots,Y_n)$ where 
  $X_1,\dots,X_n\in \alg{B}$ and  $Y_1,\dots,Y_n\in \alg{C}$ satisfies
  \begin{align*}
  \phi(P(X_1,\dots,X_n,Y_1,\dots,Y_n))
  &= \tilde{\phi}(P(X_1^{(1)},\dots,X_n^{(1)},Y_1^{(1)},\dots,Y_n^{(1)}))  \\
  &= \tilde{\phi}(P(X_1^{(1)},\dots,X_n^{(1)},Y_1^{(2)},\dots,Y_n^{(2)}))  
  \end{align*}
\end{Example}
A few remarks are in place here.
\begin{Remark}
  \begin{enumerate}
   \item 
    For classical (or free) independence,
    Definition~\ref{def:GeneralCumulants:independence} reduces to the well known fact
    that if $(X,Y)$ is a random vector with independent entries and
    $(X',Y')$ and $(X'',Y'')$ are i.i.d.\ copies, then the joint distributions
    of $(X,Y)$ and $(X',Y'')$ coincide.
   \item Note that we do not require the algebras $\alg{A}_j$ to be disjoint.
    The reader should be warned that the notion of \exch-independence is very weak
    and sometimes the term ``independence'' not even justified.
    Given a noncommutative probability space~$(\alg{A},\phi)$
    one can for instance consider the trivial exchangeability system $\alg{U}=\alg{A}$
    with the identical embedding, so that all $\alg{A}_i$ are the same and therefore
    any two subalgebras are \exch-independent.
   \item Most of the following considerations work for general multilinear maps
    into some vector space which satisfy an analogous invariance condition,
    but we did not pursue this direction yet.
   \item Contrary to the case of classical, free and boolean probabilities
    we do not have a ``free product'' construction in general, but rather
    assume that an infinite family of exchangeable subalgebras of some
    ``big'' algebra is given a priori.
   \item \exch{}-Independent algebras can be obtained in the following way. 
    Given an infinite family 
    $(\alg{A}_i)_{i\in\IN}$ of interchangeable subalgebras
    of a fixed noncommutative probability space $(\alg{U},\phi)$, 
    we fix a number~$N$ and relabel the sequence to 
    $(\alg{A}_{ij})_{i\in\IN_0,j=1,\dots, N}$.
    Let $\tilde{\alg{A}}_i=\alg{A}_{i1}\vee\alg{A}_{i2}\vee\cdots\vee\alg{A}_{iN}$ be the algebras generated by these ``clusters''
    and set $\tilde{\alg{A}}=\tilde{\alg{A}}_0$.
    Then with the embeddings
    $\iota_n:\tilde{\alg{A}}\to\tilde{\alg{A}}_n$
    we have an exchangeability system
    $\exch=(\alg{U},\phi,\alg{J})$ for $(\tilde{\alg{A}},\phi)$
    and the subalgebras $\alg{B}_j=\alg{A}_{0j}$
    are clearly \exch-exchangeable.
  \end{enumerate}
\end{Remark}

\section{Cumulants}

\subsection{Good's formula}


\begin{Definition}
  \label{def:GeneralCumulants:GoodFormula}
  Let $(\alg{A},\phi)$ be a noncommutative probability space
  and $\exch= (\alg{U},\tilde{\phi},\alg{J})$ be an exchangeability system 
  for $(\alg{A},\phi)$.
  Let $X_1,X_2,\dots,X_n\in\alg{A}$ be given random variables and
  let $X_j^{(k)}$, $k=1,\dots,n$ be their interchangeable copies.
  Let~$\omega$ be an~$n$-th primitive root of unity
  (e.g., $\omega=e^{2\pi i/n}$) and set
  \begin{equation}
    \label{eq:GeneralCumulants:Tjomega}
    X_j^{\omega} = \omega X_j^{(1)} + \omega^2 X_j^{(2)} +\dots + \omega^n X_j^{(n)}
    .
  \end{equation}
  We define the $n$th cumulant to be
  \begin{equation}
    \label{eq:generalcumulants:definition}
    K^\exchm_n(X_1, X_2, \dots, X_n)
    = \frac{1}{n}
      \,
      \tilde{\phi}(X_1^\omega X_2^\omega \dotsm X_n^\omega)
    .
  \end{equation}
\end{Definition}
The notation~``$X^{(n)}$'' and~``$X^\omega$'' instead of something like~$\Phi_n(X)$
and~$\tilde{\Phi}_n(X)$ is by no means a perfect one and may be confusing at first,
however we sticked to it because we believe that the resulting compactness of the formulae 
increases their readability. 

Next we derive the fundamental properties which justify the name ``cumulants''.
The cumulant functions are clearly multilinear.
The vanishing of ``mixed'' cumulants is almost immediate:
\begin{Proposition}
  \label{prop:GeneralCumulants:independenceimpliesmixedcumulants}
  Mixed cumulants vanish. That is, if there is a nontrivial subset $I\subseteq[n]$
  (i.e., $I\ne\emptyset$ and $I\ne[n]$) s.t.\
  $(X_j)_{j\in I}$ and $(X_j)_{j\in [n]\setminus I}$ are
  \exch{}-independent, then $K^\exchm_n(X_1,X_2,\dots,X_n)=0$.
\end{Proposition}
\begin{proof}
  \exch{}-independence implies that if we replace $(X_j^{(k)})_{j\in I}$
  by $(X_j^{(\sigma(k))})_{j\in I}$, where $\sigma\in\SG_n$ is any
  permutation, then the expectation on the right hand side
  of~\eqref{eq:generalcumulants:definition} does not change.
  This can be seen by expanding the right hand side,
  \begin{equation}
    \label{eq:GeneralCumulants:MixedCumulantExpansionProof}
    \frac{1}{n}
    \,
    \tilde{\phi}(X_1^\omega X_2^\omega \dotsm X_n^\omega)
    = \frac{1}{n} 
      \sum_{k_1,\dots,k_n=1}^n
       \tilde{\phi}(\omega^{k_1} X_1^{(k_1)} \dotsm \omega^{k_n} X_n^{(k_n)})
    ,
  \end{equation}
  and observing that to each summand
  $\phi(\omega^{k_1} X_1^{(k_1)} \dotsm \omega^{k_n} X_n^{(k_n)})$ 
  there corresponds a map $h:j\mapsto k_j$ whose kernels $\ker h|_I$ and
  $\ker h|_{[n]\setminus I}$ do not change when we apply a permutation
  $\sigma$
  to the values $\{k_j:j\in I\}$ only.
  Let's take $\sigma=(1,2,\dots,n)$ to be the full cycle.
  That is, we replace $X_j^{(k)}$ by $X_j^{(\sigma(k))}$ for
  each $j\in I$.
  Then the expectations in the sum on the right hand side of
  \eqref{eq:GeneralCumulants:MixedCumulantExpansionProof} do not change,
  but on the other hand, for $j\in I$,
  the random variable~$X_j^\omega=\sum\omega^k X_j^{(k)}$ is permuted to
  $\tilde{X}_j^\omega=\sum\omega^{k-1} X_j^{(k)}$
  and this is equal to $\bar\omega X_j^\omega$.
  Thus we can factor out $\bar\omega$ from each $X_j^{k}$ for which $j\in I$
  and get
  $$
  K^\exchm_n(X_1, X_2, \dots, X_n) = \bar\omega^{\abs{I}} K^\exchm_n(X_1, X_2, \dots, X_n)
  .
  $$
  Since we assumed $0<|I|<n$, the factor $\bar\omega^{\abs{I}}\ne1$
  and the cumulant must vanish.
\end{proof}
A converse of this proposition is also true, but for the proof the partition
lattice formulation is needed, see
Proposition~\ref{prop:GeneralCumulants:individualmixedcumulantsimpliesindependence}.

By multilinear expansion of the cumulant we immediately get the most prominent
property of cumulants, namely additivity for sums of independent variables.
\begin{Corollary}
  In the setting of Proposition~\ref{prop:GeneralCumulants:independenceimpliesmixedcumulants},
  let $(X_i)$ and $(Y_i)$ be \exch{}-independent families of noncommutative random variables.
  Then
  $$
  K^\exchm_n(X_1+Y_1,X_2+Y_2,\dots,X_n+Y_n)
  = K^\exchm_n(X_1,X_2,\dots,X_n) + K^\exchm_n(Y_1,Y_2,\dots,Y_n)
  .
  $$
\end{Corollary}

The following lemma is obvious, yet it
will turn out to be the most useful feature of Good's construction.
\begin{Lemma}
  \label{lem:GeneralCumulants:subwordexpectationvanishes}
  Let~$X_j^\omega$ be as in~\eqref{eq:GeneralCumulants:Tjomega} and
  let $j_1,j_2,\dots,j_m$ be a subsequence of $\{1,2,\dots,n\}$ with $m<n$.
  Then
  $$
  \tilde{\phi}(X_{j_1}^\omega X_{j_2}^\omega \cdots X_{j_m}^\omega)
  = 0
  $$
\end{Lemma}

\ddate{02.12.2000}

\subsection{Partition lattice formulation}
Up to now we have found the general form of properties
\eqref{eq:GeneralCumulants:axiom1} and \eqref{eq:GeneralCumulants:axiom2}.
Property \eqref{eq:GeneralCumulants:axiom3} in the form stated does not hold
in general, but in a rather weaker form which is the subject of this section.
Formula \eqref{eq:generalcumulants:definition} can be expanded and after collecting
terms we obtain the well known partition lattice formulation of the
moment-cumulant formula,
see \cite{Schutzenberger:1947:certains,Speed:1983:cumulantsI}.
\begin{Theorem}
  \label{thm:GeneralCumulants:PartitionExpansion}
  In the setting of Definition~\ref{def:GeneralCumulants:GoodFormula}
  we have
  \begin{equation}
    \label{eq:GeneralCumulants:PartitionExpansion}
    K^\exchm_n(X_1,X_2,\dots,X_n)
    = \sum_{\pi\in\Pi_n} \phi^\exchm_\pi(X_1,X_2,\dots,X_n)\,\mu(\pi,\hat{1}_n)
  \end{equation}
  where~$\mu(\pi,\sigma)$ is the M\"obius function on the partition lattice.
\end{Theorem}
\begin{proof}
\begin{align*}
  K^\exchm_n(X_1,X_2,\dots,X_n)
  &= \frac{1}{n}
     \sum_{k_1,\dots,k_n=1}^n \phi(\omega^{k_1} X_1^{(k_1)}\omega^{k_2} X_2^{(k_2)} \dotsm \omega^{k_n} X_n^{(k_n)})\\
  &= \frac{1}{n}
     \sum_{g:[n]\to [n]}
      \omega^{g(1)+g(2)+\dots+g(n)} \phi(X_1^{(g(1))} X_2^{(g(2))}\dotsm X_n^{(g(n))})\\
  &= \frac{1}{n}
     \sum_{\pi\in \Pi_n}
     \sum_{ \ker g=\pi} \omega^{g(1)+g(2)+\dots+g(n)}
     \phi^\exchm_\pi(X_1,X_2,\dots, X_n)
\end{align*}
because by assumption the value $\phi(X_1^{(g(1))}X_2^{(g(2))}\dotsm X_n^{(g(n))})$
 only depends on $\ker g$.
So we need to evaluate the function
$$
F(\pi) = \sum_{\ker g = \pi} \omega^{g(1)+g(2)+\dots+g(n)}
.
$$
To this end define another function~$G$ on $\Pi_n$ by
$$
G(\pi) = \sum_{\sigma\geq \pi} F(\sigma)
       = \sum_{\ker g\geq \pi} \omega^{g(1)+g(2)+\dots+g(n)}
.
$$
The condition $\ker g\geq \pi$ means that $g$ is constant on the blocks of
$\pi$. Since everything is commutative, for $\pi=\{\pi_1,\pi_2,\dots,\pi_p\}$
we have
$$
G(\pi) = \sum_{k_1,k_2,\dots,k_p=1}^n
          \omega^{\abs{\pi_1}k_1}
          \omega^{\abs{\pi_2}k_2} 
          \cdots
          \omega^{\abs{\pi_p}k_3} 
       = \biggl(
           \sum_{k=1}^n \omega^{\abs{\pi_1}k}
         \biggr)
         \biggl(
           \sum_{k=1}^n \omega^{\abs{\pi_2}k}
         \biggr)
         \cdots
         \biggl(
           \sum_{k=1}^n \omega^{\abs{\pi_p}k}
         \biggr)
$$
and this vanishes unless $\pi=\hat{1}_n$, in which case it equals $n$.
Therefore by Möbius inversion we obtain
$$
F(\pi) = \sum_{\sigma\geq\pi} \mu(\pi,\sigma)\,G(\sigma) = n\, \mu(\pi,\hat{1}_n)
.
$$
\end{proof}

\ddate{11.12.2000}
\subsection{Good's formula for partitioned cumulants}
The previous formula naturally leads to the definition of \emph{partitioned cumulants}.
\begin{Definition}
  Let $(\alg{A},\phi)$ be a noncommutative probability space
  and $\exch= (\alg{U},\phi,\alg{J})$ be an exchangeability system 
  for $(\alg{A},\phi)$.
  Given random variables~$X_1,X_2,\dots,X_n\in\alg{A}$
  and a partition $\pi\in\Pi_n$ we define the partitioned cumulant
  $$
  K^\exchm_\pi(X_1,X_2,\dots,X_n)
  = \sum_{\sigma\le\pi} 
     \phi^\exchm_\sigma(X_1,X_2,\dots,X_n)\,\mu(\sigma,\pi)
  .
  $$
\end{Definition}
By M\"obius inversion the general moment-cumulant formula now follows.
\begin{Proposition}
  \label{prop:GeneralCumulants:K2M}
  \begin{equation}
    \phi(X_1X_2\dotsm X_n)
    = \sum_\pi K^\exchm_\pi(X_1,X_2,\dots,X_n)
    ;
  \end{equation}
  more generally, for any partition $\sigma\in\Pi_n$ we have
  \begin{equation}
    \label{eq:generalcumulants:K2M}
    \phi^\exchm_\sigma(X_1,X_2,\dots, X_n)
    = \sum_{\pi\leq\sigma} K^\exchm_\pi(X_1,X_2,\dots,X_n)
    .
  \end{equation}
\end{Proposition}

Formula~\eqref{eq:generalcumulants:K2M} as it stands is just a reformulation
of the definition of $K^\exchm_\pi$ and M\"obius inversion.
However there is a Good type formula which can be proved
in the same way as Theorem~\ref{thm:GeneralCumulants:PartitionExpansion}.
\begin{Proposition}
  \label{prop:GeneralCumulants:PartitionGoodFormula}
  Given noncommutative random variables~$X_1$,~$X_2$, \ldots,~$X_n$
  and a partition~$\pi\in\Pi_n$,
  we choose for each $k\in [n]$ an exchangeable
  copy~$\{X_j^{(k)}:j\in [n]\}$ of the given family~$\{X_j:j\in [n]\}$,
  for each block~$B=\{k_1<k_2<\dots<k_b\}\in\pi$ we pick
  a primitive root of unity~$\omega_b$ of order~$b=|B|$, 
  and set for each~$i\in B$
  $$
  X_i^{\pi,\omega}
  = \omega_b X_i^{(k_1)}
    +
    \omega_b^2 X_i^{(k_2)}
    +
    \dots
    +
    \omega_b^{b} X_i^{(k_b)}
    ;
  $$
  that is, we do the construction of
  Definition~\ref{def:GeneralCumulants:GoodFormula}
  for each block of~$\pi$ independently.
  Then
  \begin{equation}
    \label{eq:GeneralCumulants:PartitionGoodFormula}
    K^\exchm_\pi(X_1,X_2,\dots,X_n)= \frac{1}%
                                   {\prod \abs{B}}
                              \,
                              \phi(X_1^{\pi,\omega}X_2^{\pi,\omega}\dotsm X_n^{\pi,\omega})
  \end{equation}
\end{Proposition}
Also the analogue of
Proposition~\ref{prop:GeneralCumulants:independenceimpliesmixedcumulants}
holds for partitioned cumulants.
\begin{Corollary}
  Mixed partitioned cumulants vanish:
  Let $\alg{B}$ and $\alg{C}$ be \exch-independent algebras
  and $X_i\in\alg{B}\cup\alg{C}$ some noncommutative random variables
  taken from their union.
  Let~$\pi\in\Pi_n$ be an arbitrary partition.
  If there is a block of $\pi$ which contains indices $i$ and $j$
  such that $X_i\in\alg{B}$ and $X_j\in\alg{C}$,
  then $K^\exchm_\pi(X_1,X_2,\dots,X_n)$ vanishes.
\end{Corollary}

The converse holds, too.

\begin{Proposition}
  \label{prop:GeneralCumulants:mixedcumulantsimpliesindependence}
  Let~$\alg{B}$, $\alg{C}$ be subalgebras of~$\alg{A}$ such that mixed cumulants
  vanish, that is, if for any family of random variables
  $X_1$,~$X_2$,\ldots,$X_n\in \alg{B}\cup\alg{C}$ 
  and any partition~$\pi$ such that in one of the blocks of $\pi$ there appears at least
  one element from either algebra, the cumulant
  $K^\exchm_\pi(X_1,X_2,\dots,X_n)$ vanishes,
  then $\alg{B}$ and $\alg{C}$ are \exch-independent.
\end{Proposition}
\begin{proof}
  We need to check that any finite sequence of
  random variables~$X_1,X_2,\dots,X_n\in\alg{B}\cup\alg{C}$
  satisfy the condition of Definition~\ref{def:GeneralCumulants:independence}.
  Let~$\rho$ be the partition $\{I_\alg{B},I_\alg{C}\}$
  induced by the subsets $I_\alg{B}=\{i:X_i\in\alg{B}\}$
  and~$I_\alg{C}=\{i:X_i\in\alg{C}\}$.
  Then for any partition $\pi$ we have
  $$
  \phi^\exchm_\pi(X_1,X_2,\dots,X_n)
  = \sum_{\sigma\leq\pi}
     K^\exchm_\sigma(X_1,X_2,\dots,X_n)
  = \sum_{\sigma\leq\pi\wedge\rho}
     K^\exchm_\sigma(X_1,X_2,\dots,X_n)
  $$
  because by assumption $K^\exchm_\sigma(X_1,X_2,\dots,X_n)=0$
  unless $\sigma\leq\rho$.
  Therefore
  $$
  \phi^\exchm_\pi(X_1,X_2,\dots,X_n)
  =\phi^\exchm_{\pi'}(X_1,X_2,\dots,X_n)
  $$
  whenever $\pi\wedge\rho=\pi'\wedge\rho$
  and by remark~\ref{rem:GeneralCumulants:Independence}
  the claim follows.
\end{proof}
An analogous result holds
(Proposition~\ref{prop:GeneralCumulants:individualmixedcumulantsimpliesindependence})
for concrete elements
but for the proof we need the product formula of Leonov and Shiryaev, see
Proposition~\ref{prop:GeneralCumulants:Productformula} below.

\section{Basic transformations}
\label{sec:GeneralCumulants:BasicTransformations}

In this section we investigate the behaviour of cumulants under certain
transformations of the random variables.

\subsection{Affine transformations}
\begin{Proposition}
  Let $Y_i = \sum_j \alpha_{ij} X_j + \beta_i$, $i=1,\dots,m$ be an affine transformation of
  $X_1,\dots,X_n$, then all cumulants except the
  first one do not depend on the constants $\beta_i$
  and we have for $m\geq 2$
  $$
  K^\exchm_m(Y_1,\dots,Y_m)
  = \sum_{j_1,\dots,j_m}
     \alpha_{1,j_1} \dotsm \alpha_{m,j_m}
     K^\exchm_m(X_{j_1},\dots,X_{j_m})
  $$
  An analogous formula holds for partitioned cumulants
  $K^\exchm_\pi$ if $\pi$ contains no singleton.
\end{Proposition}
\begin{proof}
  Simply expand Good's formula~\eqref{eq:generalcumulants:definition}
  multilinearly and notice that $\beta_i^{(k)}=\beta_i$ and
  consequently we have $\sum_k \omega^k \beta_i^{(k)}=0$ if $m\geq2$.
\end{proof}
\subsection{Cumulants of products}
The formula of Leonov and Shiryaev for cumulants of products 
\cite{LeonovShiryaev:1959:method}
and Speed's proof 
\cite[Proposition~4.3]{Speed:1983:cumulantsI}%
\ppaper[Prop.~4.3]{Speed}{Cumulants and partition lattices}{Austral. J. Statist. 25 (1983) 378--388}
can immediately be transferred to the noncommutative case.
The analogous formula for free cumulants was found with a different proof
and many applications in \cite{KrawczykSpeicher:2000:combinatorics},
see also \cite{Speicher:2000:conceptual,Cabanal-Duvillard:1999:noncrossing} for other proofs.
\begin{Definition}
  \label{def:generalcumulants:indecomposable}
  Let $(X_{i,j})_{i\in[m],j\in [n_i]}\subseteq\alg{A}$
  be a family of noncommutative random variables,
  in total $n=n_1+n_2+\dots+n_m$ variables.
  Then every partition $\pi\in\Pi_m$ induces a partition $\tilde{\pi}$
  on $[n]=\{(i,j):i\in [m], j\in[n_i]\}$ with blocks
  $\tilde{B}=\{(i,j):i\in B, j\in[n_i]\}$.
  A partition $\sigma\in\Pi_n$ is called \emph{decomposable} relative to $\pi$
  if $\sigma\le \tilde\pi$ and it is \emph{indecomposable} if no such relation holds
  other than $\sigma\le\tilde{\hat{1}}_m$;
  in other words, if $\sigma\vee\tilde{\hat{0}}_m = \tilde{\hat{1}}_m$.
\end{Definition}

\begin{Proposition}
  \label{prop:GeneralCumulants:Productformula}
  With the settings of Definition~\ref{def:generalcumulants:indecomposable}
  we have
  $$
  K^\exchm_m(\prod_{j_1} X_{1,j_1},\prod_{j_2} X_{2,j_2},\dots,\prod_{j_m} X_{m,j_m})
  =\sum_{\sigma\in\Pi_n\ \text{indec.}}
    K^\exchm_\sigma(X_{1,1},X_{1,2},\dots,X_{1,n_1},X_{2,1}\dots,X_{m,n_m})
  $$
\end{Proposition}
\begin{proof}
  Denote $\tilde{X}_i = \prod_j X_{i,j}$.
  Then using the fact that mixed cumulants vanish we have for $\pi\in\Pi_m$ that
  $$
  F(\pi)
  :=\phi^\exchm_\pi(\tilde{X}_1,\tilde{X}_2,\dots,\tilde{X}_m)
  = \sum_{\sigma\le\tilde\pi} K^\exchm_\sigma(X_1,X_2,\dots,X_n)
  $$
  Define for $\pi\in\Pi_m$ the partial sum
  $$
  f(\pi) = \sum_{\substack{\sigma\le\tilde\pi\\ \sigma\not\le \tilde\rho\, \forall \rho<\pi}}
            K^\exchm_\sigma(X_1,X_2,\dots,X_n)
         = \sum_{\substack{
                   \sigma\le\tilde\pi\\
                   \sigma\vee\tilde{\hat0}_m= \tilde\pi}}
            K^\exchm_\sigma(X_1,X_2,\dots,X_n)
  $$
  Then obviously $F(\pi)=\sum_{\rho\le\pi} f(\rho)$ and by M\"obius inversion
  $f(\pi)=\sum_{\sigma\le\pi} F(\sigma)\,\mu(\sigma,\pi)$.
  Therefore
  \begin{align*}
    K^\exchm_m(\tilde{X}_1,\tilde{X}_2,\dots,\tilde{X}_m)
    &= \sum_\pi \phi^\exchm_\pi(\tilde{X}_1,\tilde{X}_2,\dots,\tilde{X}_m)\,\mu(\pi,\hat{1}_m) \\
    &= \sum_\pi F(\pi)\,\mu(\pi,\hat{1}_m)\\
    &= f(\hat{1}_m)
     = \sum_{\sigma\ \text{indec.}}
        K^\exchm_\sigma(X_{1,1},X_{1,2},\dots,X_{m,n_m})
  \end{align*}
\end{proof}

An analogous formula holds for partitioned cumulants.
\begin{Proposition}
  With the settings of Definition~\ref{def:generalcumulants:indecomposable}
  we have
  $$
  K^\exchm_\pi(\prod_{j_1} X_{1,j_1},\prod_{j_2} X_{2,j_2},\dots,\prod_{j_m} X_{m,j_m})
  =\sum_{\substack{\sigma\in\Pi_n\\ \sigma\vee \tilde{\hat0}_m=\tilde\pi}}
    K^\exchm_\sigma(X_{1,1},X_{1,2},\dots,X_{m,n_m})
  $$
\end{Proposition}
Now we are able to prove a stronger version of 
Proposition~\ref{prop:GeneralCumulants:mixedcumulantsimpliesindependence}
\begin{Proposition}
  \label{prop:GeneralCumulants:individualmixedcumulantsimpliesindependence}
  A family of random variables $X_1,X_2,\dots,X_m\in\alg{A}$ is \exch-independent
  if and only if mixed cumulants vanish, i.e.,
  if for every finite sequence $X_{i_1},X_{i_2},\dots,X_{i_n}$ taken
  from the family and for every partition $\pi\in\Pi_n$, such that
  some block of $\pi$ contains two different indices, the cumulant
  $K^\exchm_\pi(X_{i_1},X_{i_2},\dots,X_{i_n})$ vanishes.
\end{Proposition}
\begin{proof}
  This can be reduced to
  Proposition~\ref{prop:GeneralCumulants:mixedcumulantsimpliesindependence},
  with the help of the product formula.
  For simplicity we consider the case of two random variables $X_1$ and $X_2$ only.
  All that is left to show is that if all mixed cumulants of $X_1$ and $X_2$ vanish,
  then all mixed cumulants of elements from the algebras $\alg{B}_1$
  and $\alg{B}_2$ generated by $X_1$ and $X_2$, respectively, vanish.
  In other words, we have to show that mixed cumulants of polynomials 
  $P_1(X_{i_1}), P_1(X_{i_2}),\dots,P_1(X_{i_m})$ vanish.
  By multilinearity it suffices to consider simple powers
  $X_{i_1}^{k_1},X_{i_2}^{k_2},\dots,X_{i_m}^{k_m}$.
  Let $n=k_1+k_2+\dots+k_m$ and denote $\tilde{\hat0}_m$ the partition induced
  on $[n]$ (see Definition~\ref{def:generalcumulants:indecomposable}).
  In other words, $\tilde{\hat0}_m$ is the partition with interval blocks 
  $I_1$, $I_2$,\ldots, $I_m$ of length $k_1$, $k_2$,\ldots, $k_m$.
  Let $\pi\in\Pi_m$ and assume that some block of $\pi$ contains elements $X_1^{k_r}$
  and $X_2^{k_s}$.
  By Proposition~\ref{prop:GeneralCumulants:mixedcumulantsimpliesindependence}
  we have
  \begin{multline*}
  K^\exchm_\pi(X_{i_1}^{k_1},X_{i_2}^{k_2},\dots,X_{i_m}^{k_m})\\
  =\sum_{\substack{\sigma\in\Pi_n\\ \sigma\vee \tilde{\hat0}_m=\tilde\pi}}
    K^\exchm_\sigma(X_{i_1},X_{i_1},\dots,X_{i_1},
             X_{i_2},X_{i_2},\dots,X_{i_2},
             \dots,
             X_{i_m},X_{i_m},\dots,X_{i_m})
  \end{multline*}
  By assumption the blocks $I_r$ and $I_s$ are contained in one block of $\tilde\pi$ and
  therefore a partition $\sigma\in\Pi_n$ which satisfies $\sigma\vee\tilde{\hat0}_m=\tilde\pi$
  must connect at least one element from each $I_r$ and $I_s$.
  This implies that some block of~$\sigma$ must contain both~$X_1$ and~$X_2$,
  that is, it is a mixed cumulant of $X_1$ and $X_2$, which by assumption vanishes.
\end{proof}

\subsection{Cumulants of matrices}
\label{sec:GeneralCumulants:matrices}
\ddate{10.04.2001}
Free cumulants of matrices with free entries are computed in
\cite{NicaShlyakhtenkoSpeicher:2001:Rcyclic}.
\ppaper[sec.~6]{Nica/Shlyakhtenko/Speicher}{R-cyclic families of matrices in free probability}{MSRI 2000-038, arXiv:math.OA/0101025}
Using Good's formula it is actually quite simple to obtain a formula for the
cumulants of matrices of random variables.

\begin{Proposition}
  Let $(\alg{A},\phi)$ be a noncommutative probability space
  and $\exchm=(\alg{U},\phi,\alg{J})$ be an exchangeability system
  for $(\alg{A},\phi)$.
  Then $M_d(\exchm)=(M_d(\alg{U}),\psi,M_d(\alg{J}))$
  is an exchangeability system for matrix-valued probability space
  $(M_d(\alg{A}),\psi)$
  where $\psi= Id_{M_d}\ox\phi$
  is the conditional expectation 
  from $M_d(\alg{A})$
  onto be the subalgebra of constant matrices
  $\IM_d\subseteq M_d(\alg{A})$ 
  given by~$\psi([X_{i,j}]) = [\phi(X_{i,j})]$.
  The subalgebras $M_d(\alg{A}_i)$ are clearly interchangeable
  with respect to $\psi$
  if the algebras $\alg{A}_i$ are interchangeable with respect to
  $\phi$
  and the cumulants of the matrices~$X_k = [X_{i,j}(k)]\in M_d(\alg{A}_0)$
  are matrices with entries
  $$
  K^{M_d(\exchm)}_\pi(X_1,X_2,\dots,X_n)_{i,j}
  =\sum_{i_1,\dots,i_{n-1}=1}^d
    K^\exchm_\pi(X_{i,i_1}(1),X_{i_1,i_2}(2),\dots,X_{i_{n-1},j}(n))
  $$
\end{Proposition}
Theorem~6.2 of \cite{NicaShlyakhtenkoSpeicher:2001:Rcyclic} is a consequence
of this because 
$M_d(\alg{A}_i)$ are free with amalgamation over $\IM_d$
if and only if $\alg{A}_i$ are free.

\ddate{11.04.2001}
Cumulants of matrices are a special case of cumulants of tensor products.
The above observation is a special case of the following proposition.
\begin{Proposition}
  Let $\alg{A}_i\subseteq\alg{U}$ be interchangeable 
  with respect to $\phi$ and
  let $\alg{C}$ be another algebra, then the algebras $\alg{C}\ox\alg{A}_i$
  are interchangeable w.r.\ to $\psi=E_\alg{C}=Id_\alg{C}\ox\phi$.
  and we have a new exchangeability system
  $\alg{C}\ox \exchm](\alg{C}\ox\alg{U},\psi,Id_\alg{C}\ox\alg{J})$.
  For $X_i = \sum_j C_{ij}\ox T_{ij}$ we have 
  $$
  E_{\alg{C}} (X_1\cdots X_n)
   = \sum_{j_1,j_2,\dots,j_n}
      C_{1j_1} C_{2j_2}\cdots C_{nj_n}
      \,
      \phi(T_{1j_1} T_{2j_2}\cdots T_{nj_n})
  $$
  and therefore the corresponding cumulants satisfy
  $$
  K^{\alg{C}\ox\exchm}_\pi (X_1,X_2,\dots, X_n)
   = \sum_{j_1,j_2,\dots,j_n}
      C_{1j_1} C_{2j_2}\cdots C_{nj_n}
      \,
      K^\exchm_\pi(T_{1j_1}, T_{2j_2},\dots, T_{nj_n})
  $$
\end{Proposition}
Another possibility is to choose a state~$\rho$ on $\alg{C}$
and consider the product state $\rho\ox\phi$ on $\alg{C}\ox\alg{A}$.
The corresponding cumulants are then given by
$$
K^{\alg{C}\ox\exchm,\rho\ox\phi}_\pi (X_1,X_2,\dots, X_n)
=\rho(K^{\alg{C}\ox\exchm}_\pi (X_1,X_2,\dots, X_n))
$$
In the case of matrices where~$\alg{C}=\IM_d$ a natural choice for~$\rho$
is the trace~$\tau_d=\frac{1}{d}\tr$
and in this case the cumulants are given by ``cyclic sums'' of the original cumulants:
$$
K^{M_d(\exchm),\tau_d\ox\phi}_\pi(X_1,X_2,\dots,X_n)
=\frac{1}{d}\sum_{i_1,\dots,i_{n}=1}^d
  K^\exchm_\pi(X_{i_n,i_1}(1),X_{i_1,i_2}(2),\dots,X_{i_{n-1},i_n}(n))
.
$$
The reader should be warned that for example in the case of freeness
the cumulants $K^{M_d(\exchF),\tau_d\ox\phi}_\pi(X_1,X_2,\dots,X_n)$
are \emph{different} from the free cumulants 
$K^{\exchF,\tau_d\ox\phi}_\pi(X_1,X_2,\dots,X_n)$, which
are related to a different exchangeability system and it is
much more difficult to express these in terms of the amalgamated cumulants.

The following lemma allows us to remove the identity element from cumulants.
It is an easy consequence of Proposition~\ref{prop:GeneralCumulants:PartitionGoodFormula}.
\begin{Lemma}
  \label{lem:GeneralCumulants:RemoveIdentity}
  Let $X_i\in\alg{A}$ s.t.\ $X_j=I$ for $j\in I\subseteq \{1,2,\dots,n\}$.
  Let $\pi\in\Pi_n$ have singletons at each $j\in I$ (and possibly more).
  Let $\tilde{\pi}$ be the partition obtained by removing these singletons from $\pi$.
  Then
  $$
  K^\exchm_\pi(X_1,X_2,\dots,X_n) = K^\exchm_{\tilde{\pi}}(X_i:i\in\{1,2,\dots, n\}\setminus I)
  $$
\end{Lemma}
On the other hand it is clear from Good's formula that if~$I$ appears in a block of size at least two then the corresponding cumulant~$K_\pi$ vanishes.

\subsection{The recursion formula}

\ddate{18.05.2001}
The general form of the recursion formula
\noteQ{\eqref{eq:ClassicalCumulants:Recursion}}{\eqref{eq:GeneralCumulants:ClassicalCumulantRecursion}}
(see also \cite{Speicher:1994:multiplicative} for the free analog)
is as follows. The general philosophy is to replace cumulants
by expectations of $X^\omega$.
\begin{Proposition}
  \begin{equation}
    \label{eq:generalcumulants:recursion}
  \phi(X_1 X_2\cdots X_n)
  = \sum_{\substack{A\subseteq [n]\\ A\ni 1}}
     \frac{1}{\abs{A}}\,
     \phi(X_1^{A,\omega} X_2^{A,\omega} \cdots X_n^{A,\omega} )
  \end{equation}
where $\omega$ is a root of unity of order $\abs{A}$ and
$$
X_j^{A,\omega}
= \begin{cases}
    X_j^{\omega}
    = \omega X_j^{(1)} + \omega^2 X_j^{(2)} + \dots + \omega^{\abs{A}} X_j^{(\abs{A})}
                & j\in A\\
    X_j^{(0)}   & \text{otherwise}
  \end{cases}
$$
\end{Proposition}

\begin{proof}
We have (cf.\ \eqref{eq:generalcumulants:K2M})
\begin{align*}
  \phi(X_1X_2\dotsm X_n)
  &= \sum_{\pi\in\Pi_n}
      K^\exchm_\pi(X_1,X_2,\dots,X_n)\\
  &= \sum_{\substack{A\subseteq [n]\\ A\ni 1}}
      \sum_{\pi\in\Pi_{[n]\setminus A}}
       K^\exchm_{\pi\cup\{A\}}(X_1,X_2,\dots,X_n)
\end{align*}
Now for fixed $A$ define a new multilinear functional $\phi_A$
for tuples of random variables which are indexed by $[n]\setminus A$, the complement
of $A$, namely for any such tuple
$$
\phi_A((S_j)_{j\in [n]\setminus A})
= \frac{1}{\abs{A}}\,
  \phi(S_1S_2\cdots S_n)
$$
where we fill up the sequence to an $n$-tuple by setting
$S_j = X_j^\omega$ for $j\in A$, where $(X_j^{(k)})_{j\in A}$ and $(S_j)_{j\in [n]\setminus A}$
are chosen independent.
Then exchangeable families of $(S_j)$ remain exchangeable and we can define cumulants
for this functional:
\begin{align*}
  K^{\exchm,\phi_A}_\pi((S_j)_{j\in [n]\setminus A})
  &= \frac{1}{\prod_{B\in\pi}\abs{B}}
     \,
     \phi_A((S_j^{\omega,\pi})_{j\in[n]\setminus A})\\
  &= \frac{1}{\abs{A}\prod_{B\in\pi}\abs{B}}
     \,
     \phi(\tilde S_1^{\omega,\pi} \tilde S_2^{\omega,\pi}\cdots \tilde S_n^{\omega,\pi}) \\
  &= K^\exchm_{\pi\cup\{A\}}(\tilde S_1,\tilde S_2,\dots,\tilde S_n)
\end{align*}
where 
$$
\tilde S_j
= \begin{cases}
    X_j & j \in A\\
    S_j & j \in A
\end{cases}
;
$$
now \eqref{eq:generalcumulants:recursion} follows.
\end{proof}

\subsection{Pyramidal independence}
\label{ssec:GeneralCumulants:PyramidalIndependence}

Pyramidal independence implies that cumulants are multiplicative on noncrossing partitions
and, more generally, along the connected components.
\begin{Definition}[{\cite{BozejkoSpeicher:1996:interpolations}}]
  \label{def:GeneralCumulants:PyramidalIndependence}
  Two subalgebras $\alg{B}$ and $\alg{C}$ of $\alg{A}$ 
  satisfy \emph{pyramidal independence} if $\phi(X Y X') = \phi(XX')\,\phi(Y)$
  whenever $X,X'\in\alg{B}$ and $Y\in\alg{C}$ and vice versa.

  We will say that an interchangeable family of algebras $\alg{A}_i$ satisfies
  pyramidal independence if for any choice of disjoint index sets~$I$ and~$J$
  the algebras $\alg{A}_I$ and $\alg{A}_J$, generated by $(A_i)_{i\in I}$ and 
  $(A_i)_{i\in J}$ respectively, satisfy pyramidal independence.
\end{Definition}

\begin{Proposition}
  \label{prop:GeneralCumulants:PyramidalIndependence}
  If the algebras $\alg{A}_i$ satisfy pyramidal independence, then
  the moments and consequently the \exch-cumulants factorize along the connected components.
  That is, if $\pi$ has connected components $\pi'$, $\pi''$,\ldots{} etc.,
  then
  $$
  \phi^\exchm_\pi(X_1,X_2,\dots,X_n)
  = \phi^\exchm_{\pi'}(X_1,X_2,\dots,X_n)\,
    \phi^\exchm_{\pi''}(X_1,X_2,\dots,X_n)
    \cdots
  $$
  and
  $$
  K^\exchm_\pi(X_1,X_2,\dots,X_n)
  = K^\exchm_{\pi'}(X_1,X_2,\dots,X_n)\,
    K^\exchm_{\pi''}(X_1,X_2,\dots,X_n)
    \cdots
  $$
\end{Proposition}

\subsection{Relations between different cumulants}
\label{sec:GeneralCumulants:relations}
Considering matrices (section~\ref{sec:GeneralCumulants:matrices}) and
classical exchangeable random variables
(section~\ref{sec:GeneralCumulantExamples:classicalexchangeable})
one might wonder what are the relations between the different kinds of
cumulants: given different exchangeability systems for a fixed
probability space $(\alg{A},\phi)$, is it possible to express one kind of
cumulants in terms of the other?

One can indeed express free cumulants in terms of classical cumulants, namely
\cite{Lehner:2002:connected}
$$
K_n^\exchF(X) = \sum_{\pi\in\Pi_n^{conn}} \kappa_\pi(X)
$$
where the sum runs over all connected partitions.
One can show by induction that an inverse formula holds as well, 
but there is apparently no way to write it down explicitly.

In general one cannot expect to be able to express one kind of cumulants in
terms of another. 
For example, the $q$-cumulants of some a noncommutative random variable $X$
(cf.~\cite{Lehner:2002:Cumulants3}) are not determined by the moments of $X$
alone \cite{vanLeeuwenMaassen:1996:obstruction}, but depend on the concrete
realization of $X$ as an operator on $q$-Fock space;
free or classical cumulants however only depend on the moments of $X$.
Therefore it is not possible to express $q$-cumulants in terms of free
cumulants. The converse, however, is true, because the $q$-cumulants determine
the moments of $X$ and the moments determine both free and classical cumulants.

Another question is the following.
Assume that we are given an operator-valued exchangeability system
$(\alg{E},\psi,\alg{J})$ for the operator-valued noncommutative probability space
$(\alg{A},\psi)$ with values in some subalgebra $\alg{B}$.
Choosing an arbitrary state $\phi$ on $\alg{B}$,
$(\alg{E},\phi\circ\psi,\alg{J})$ becomes an exchangeability system for
$(\alg{A},\phi\circ\psi)$ and trivially
$$
K^{\exchm,\phi}_\pi(X_1,X_2,\dots,X_n)
= \phi(K^{\exchm,\psi}_\pi(X_1,X_2,\dots,X_n))
$$
as already observed in section~\ref{sec:GeneralCumulants:matrices}, where
$\alg{B}=\IM_n$.
More interesting is the question, how to express for example
free cumulants of matrices $A_k=[a_{i,j}(k)]$ w.r.\ to $\tau_n\ox\phi$ in terms
of the free cumulants of the entries $a_{i,j}$. These are different from the
cumulants above, because freeness with amalgamation w.r.\ to $\psi$ does
not imply freeness w.r.\ to $\phi$. Some aspects of this question are
treated in \cite{NicaShlyakhtenkoSpeicher:2001:Rcyclic}.

For classical cumulants there is Brillinger's formula~\eqref{eq:Brillinger},
which expresses classical cumulants in terms of conditional cumulants.
There is a certain free analog~\cite{Lehner:2004:cumulants4},
but we were not able to find a formulation of Brillinger's formula in the general
context.



\section{Examples}
\label{sec:GeneralCumulants:Examples}

In this section we review some known facts about various cumulants in
the light of Good's formula. It is easily checked that all the examples considered
here satisfy the axioms of Definition~\ref{def:GeneralCumulants:independence}.
We start with the simplest cases, namely classical independent random
variables and conditionally independent random variables.
De Finetti's theorem (see Theorem~\ref{thm:GeneralCumulants:DeFinettiTheorem} below)
tells us that we cannot expect more examples from commutative probability theory.
Then various notions of cumulants from truly noncommutative probability spaces
are reviewed, like free, boolean, conditionally free etc.
Considerations on Fock spaces are postponed to a separate
paper~\cite{Lehner:2002:Cumulants3}.

\subsection{Classical cumulants}
Given a classical probability space $(\Omega,\Sigma,\mu)$,
we construct the noncommutative probability space $L^\infty(\Omega,\mu)$
which is commutative in this case.
The expectation is denoted as usual by $\IE$.
We can construct infinitely many interchangeable copies of 
$\alg{A}=L^\infty(\Omega,\mu)$
by embedding it 
into~$\alg{U}=L^\infty(\Omega^\infty,\mu^{\otimes\infty})$.
This gives rise to an exchangeability system for $\alg{A}$
and independence of subalgebras of $\alg{A}$ is equivalent to
exchangeability with respect to this exchangeability system.
From the very definition of classical independence it follows immediately
that for a partition $\pi\in\Pi_n$, the partitioned moment is
$$
\IE_\pi(X_1,X_2,\dots,X_n) = \prod_{B\in\pi} \IE\prod_{i\in B} X_i
$$
and similarly for the cumulants we have
$$
\kappa_\pi(X_1,X_2,\dots,X_n) = \prod_{B\in\pi} \kappa_{\abs{B}}( X_i : i \in B )
$$
and we deduce from
\eqref{eq:GeneralCumulants:PartitionExpansion} the well known formula of
Sch\"utzenberger~\cite{Schutzenberger:1947:certains}
$$
\kappa_n(X_1,X_2,\dots,X_n) = \sum_{\pi\in\Pi_n} \kappa_\pi(X_1,X_2,\dots,X_n) \, \mu(\pi,\hat1_n)
$$

\subsection{Classical conditionally independent random variables}

If we take conditional expectations to a
$\Sigma$-subalgebra~$\alg{B}\subseteq\alg{A}$,
then the partitioned conditional expectation factors
just like the partitioned expectation of independent random variables,
and the result is a $\alg{B}$-measurable random variable:
$$
m_\pi(X_1,X_2,\dots,X_n | \alg{B})
= \prod_j \IE (\prod_{i\in \pi_j} X_i | \alg{B} )
;
$$
consequently the partitioned $\alg{B}$-valued \emph{conditioned cumulants}
factorize as well and can be expressed via M\"obius inversion:
$$
\kappa_\pi(X_1,X_2,\dots,X_n | \alg{B})
 = \sum_{\sigma\leq\pi}
    m_\sigma(X_1,X_2,\dots,X_n | \alg{B})
    \,\mu(\pi,\hat1_n)
$$
Conditioned cumulants can be used to
detect conditional independence, namely if $X_1,X_2,\dots,X_n$ can be divided
into two groups which are independent conditionally on $\alg{B}$, then the 
cumulant~$\kappa_n(X_1,X_2,\dots,X_n|\alg{B})$ vanishes.

\subsection{Classical exchangeable random variables}
\label{sec:GeneralCumulantExamples:classicalexchangeable}
Classical (infinite) sequences of exchangeable random variables
are characterized by \emph{de Finetti's theorem}.
\begin{Theorem}[{De Finetti \cite{ChowTeicher:1978:probability,Kingman:1978:uses}}]
  \label{thm:GeneralCumulants:DeFinettiTheorem}
  Let $(X_i)$ be an infinite exchangeable sequence of random variables.
  Then 
  $X_i$ are i.i.d.\ 
  conditionally on some $\sigma$-algebra $\alg{B}$.  
\end{Theorem}
Note that no such theorem holds for finite sequences.

There are two different kinds of cumulants for classical exchangeable random
variables which one may consider,
namely the classical cumulants $\kappa_n(X)$ and the cumulants 
$K^\exchm_n(X)$
induced by the exchangeability relation.
By De~Finetti's Theorem, the exchangeability system can be
realized by considering conditionally independent copies
and therefore the corresponding cumulants can be expressed
in terms of the conditional cumulants $\kappa_n(X|\alg{B})$.
The conditional cumulants are
$$
\kappa_\pi(X_1,X_2,\dots,X_n|\alg{B})
= \sum_{\sigma\leq\pi}
   \IE(X_1^{(\sigma(1))} X_2^{(\sigma(2))}\dotsm X_n^{(\sigma(n))}|\alg{B})
   \,
   \mu(\sigma,\pi)
$$
and also
\begin{align*}
K^\exchm_\pi(X_1,X_2,\dots,X_n)
&= \sum_{\sigma\leq\pi}
   \IE(X_1^{\sigma(1)} X_2^{\sigma(2)} \dotsm X_n^{\sigma(n)})
   \,
   \mu(\sigma,\pi)\\
&= \IE\sum_{\sigma\leq\pi}
   \IE(X_1^{\sigma(1)} X_2^{\sigma(2)} \dotsm X_n^{\sigma(n)}|\alg{B})
   \,
   \mu(\sigma,\pi)\\
&= \IE(\kappa_\pi(X|\alg{B}))
.
\end{align*}
Note that $K^\exchm_\pi$ does not factorize along the blocks in this case.
On the other hand, the classical cumulants $\kappa_n(X)$ are
given by the more complicated formula of Brillinger
\cite{Brillinger:1969:calculation}:
\begin{equation}
  \label{eq:Brillinger}
  \kappa_n(X_1,X_2,\dots,X_n)
  = \sum_{\pi\in\Pi_n}
     \kappa_{\abs{\pi}}( \kappa_{\abs{\pi_j}}(X_i : i\in \pi_j | \alg{B})
                                        : j=1,\dots,\abs{\pi})
  .
\end{equation}

\subsection{Free cumulants}
\label{sec:GeneralCumulantExamples:free}
Free independence is one of the most fundamental notions of independence
in noncommutative probability. It was introduced by Voiculescu in
\cite{Voiculescu:1985:symmetries}, where among many other facts
existence of cumulants was shown. 
A systematic theory was established by Speicher
\cite{Speicher:1994:multiplicative},
who found the fundamental connection to the lattice of noncrossing partitions.
For further information on free probability we refer to
\cite{VDN:1992:free,
      Voiculescu:2000:lectures,
      HiaiPetz:2000:semicircle,
      NicaSpeicher:2000:combinatorics}.
Here we rederive the basic facts in an elementary way.

Let us recall the definition of free independence.
\begin{Definition}[{\cite{Voiculescu:1985:symmetries}}]
  \label{def:GeneralCumulants:Freeness}
  Let $(\alg{A},\phi)$ be a noncommutative probability space.
  Subalgebras $\alg{A}_i\subseteq\alg{A}$ are called \emph{free} if
  $\phi(X_1X_2\dotsm X_n)=0$
  whenever $\phi(X_j)=0$, $X_j\in \alg{A}_{i_j}$ and $i_j\ne i_{j+1}$.
  Elements $X_i\in\alg{A}$ are said to be free if the algebras they generate
  are free.
\end{Definition}
It is not difficult to show that the mixed moments of free random variables depend only
on the moments of the individual random variables in a universal way.
An exchangeability system~$\alg{F}$
can be constructed by taking the reduced free product of copies of a given algebra:
Let $(\alg{A},\phi)$ be a noncommutative probability space,
$\alg{U}=\bigstar_i \alg{A}_i$
be the free product of infinitely many copies of $\alg{A}$
and $\tilde\phi = \bigstar_i \phi_i$ the free product state
\cite{Voiculescu:1985:symmetries,Avitzour:1982:free}.
Then the~$\alg{A}_i$ are interchangeable copies of~$\alg{A}$
and two subalgebras $\alg{B},\alg{C}\subseteq\alg{A}$ are free
if and only if they are $\exchF$-independent in the sense of
Definition~\ref{def:GeneralCumulants:independence}.

It is easily seen from Definition~\ref{def:GeneralCumulants:Freeness} that
freeness implies pyramidal independence 
(see section~\ref{ssec:GeneralCumulants:PyramidalIndependence})
and therefore by Proposition~\ref{prop:GeneralCumulants:PyramidalIndependence}
we have factorization along noncrossing partitions.

\begin{Proposition}
  \label{prop:GeneralCumulants:FreeNoncrossingExpectationFactors}
  For a noncrossing partition $\pi=\{\pi_1,\pi_2,\dots,\pi_p\}\in\NC_n$
  the partitioned expectations and cumulants factorize:
  \begin{align*}
    \phi^\exchF_\pi(X_1,X_2,\dots,X_n)
    &= \prod_{B\in\pi} \phi(\arrowprod_{i\in B} X_i)\\
    K^\exchF_\pi(X_1,X_2,\dots,X_n)
    &= \prod_{B\in\pi} K^\exchF_{\abs{B}}( X_i : i\in B)
  \end{align*}
  The products and sequences are to be taken in the order of the indices.
\end{Proposition}

The expression for $\phi^\exchF_\pi$ is rather complicated if $\pi$ has a crossing.
The cumulants however vanish in this case.
\begin{Proposition}
  \label{prop:GeneralCumulantExamples:free:crossingcumvanish}
  If $\pi\in\Pi_n$ has a crossing then
  $$
  K^\exchF_\pi(X_1,X_2,\dots,X_n)=0
  $$
  for any choice of $X_i$.
\end{Proposition}
\begin{proof}
  We use Proposition~\ref{prop:GeneralCumulants:PartitionGoodFormula}.
  By pyramidal independence
  we can factor out the connected components of $\pi$.
  As $\pi$ has a crossing, there is at least one connected component
  which is not a block itself, i.e., it contains at least $2$ blocks.
  It is enough to show that the contribution of this connected component is
  zero. 
  So without loss of generality we may assume that $\pi$ is connected.
  In this case no block of $\pi$ is an interval
  because of Lemma~\ref{lem:GeneralCumulants:subwordexpectationvanishes}
  we find that the cumulant~$K^\exchF_\pi(X_1,X_2,\dots,X_n)$ equals the expectation
  of an alternating word of centered free random variables.
  Therefore it vanishes.
\end{proof}
It follows that
$$
\phi^\exchF_\pi(X_1,X_2,\dots,X_n)
= \sum_{\substack{\sigma\leq\pi\\ \sigma\in\NC_n}}
   K^\exchF_\sigma(X_1,X_2,\dots,X_n)
$$
and we can apply M\"obius inversion on the lattice of noncrossing
partitions to obtain Speicher's formula~\cite{Speicher:1994:multiplicative}:
For a noncrossing partition~$\pi$ we have
$$
K^\exchF_\pi(X_1,X_2,\dots,X_n)
= \sum_{\substack{\sigma\in\NC_n\\ \sigma\leq\pi}}
   \phi^\exchF_\sigma(X_1,X_2,\dots,X_n)\, \mu_{\NC}(\sigma,\pi)
$$
where $\mu_\NC(\pi,\sigma)$ is the M\"obius function on the lattice of
noncrossing partitions, which was found by
Kreweras~\cite{Kreweras:1972:partitions}, 
and~$\phi^\exchF_\sigma$ factorizes according to
Proposition~\ref{prop:GeneralCumulants:FreeNoncrossingExpectationFactors}.

\begin{Remark}
  It follows from the considerations above that we have the 
  remarkable identity
  $$
  \sum_{\pi\in\NC_n}
   \phi^\exchF_\pi(X_1,X_2,\dots,X_n)\, \mu_{\NC}(\pi,\hat1_n)
  =
  \sum_{\pi\in\Pi_n}
   \phi^\exchF_\pi(X_1,X_2,\dots,X_n)\, \mu(\pi,\hat1_n)
  $$
  (both sides are equal to $K^\exchF_n(X_1,X_2,\dots,X_n)$),
  for which there is probably no simple ``direct'' proof.
\end{Remark}

\subsection{Operator valued free cumulants}
\label{sec:GeneralCumulants1:Examples:Amalgamated}
There is an operator valued generalization of free probability
which was also developed by Voiculescu 
\cite{Voiculescu:1985:symmetries,Voiculescu:1995:operations}.
Roughly speaking, operator valued free probability is obtained
by replacing the field $\IC$ by a subalgebra $\alg{B}$ of the given
algebra and the expectation map $\psi$, which has values in $\alg{B}$,
can be seen as an analogue of conditional expectations in classical
probability.
\begin{Definition}[{\cite{Voiculescu:1985:symmetries}}]
  Let $(\alg{A},\psi)$ be a $\alg{B}$-valued
  noncommutative probability space, that is, $\alg{B}$ is a unital
  subalgebra of $\alg{A}$ and $\psi:\alg{A}\to\alg{B}$ is a 
  conditional expectation.
  Subalgebras $\alg{A}_i\subseteq\alg{A}$ are called 
  \emph{free (with amalgamation) over $\alg{B}$} or \emph{$\alg{B}$-free}
  if $\psi(X_1X_2\dotsm X_n)=0$
  whenever $\psi(X_j)=0$, $X_j\in \alg{A}_{i_j}$ and $i_j\ne i_{j+1}$.
  Elements $X_i\in\alg{A}$ are said to be $\alg{B}$-free if the algebras they
  generate are $\alg{B}$-free.
\end{Definition}
An exchangeability system~$\exchFA$ realizing freeness with amalgamation
can be constructed by taking
amalgamated free products of algebras.
The corresponding cumulants are again governed by the lattice
of noncrossing partitions as found by Speicher
\cite{Speicher:1998:combinatorial}.
A ``nested'' analogue of pyramidal independence holds and
by a similar argument as above we have a factorization of 
partitioned expectations along connected components;
this time, however, the factors are noncommutative and remain nested.
\begin{Proposition}
  Let $\pi\in\Pi_n$ be an arbitrary partition and let
  $\sigma$ be a connected component of $\pi$,
  such that $\bigcup \sigma_j=\{k,k+1,\dots,l\}$ is an interval.
  Then
  $$
  \psi_\pi(X_1,X_2,\dots,X_n)
  = \psi_{\pi\setminus\sigma}(X_1,
                             X_2,
                             \dots,
                             X_{k-1},
                             \psi_\sigma(X_k X_{k+1}\dotsm X_l)X_{l+1},
                             X_n)
  $$
  In particular, if
  $\pi=\{\pi_1,\pi_2,\dots,\pi_n\}\in\NC_n$ is a noncrossing partition
  and $\pi_j=\{k,k+1,\dots,l\}\in\pi$ is an interval block,
  then 
  $$
  \psi_\pi(X_1,X_2,\dots,X_n)
  = \psi_{\pi\setminus\pi_j}(X_1,
                             X_2,
                             \dots,
                             X_{k-1},
                             \psi(X_k X_{k+1}\dotsm X_l)X_{l+1},
                             X_n)
  $$
\end{Proposition}

The noncrossing cumulants enjoy the same factorization property,
while the crossing cumulants vanish.
The proof is essentially the same as above.
\begin{Proposition}
  If $\pi\in\Pi_n$ has a crossing then
  $K^\exchFA_\pi(X_1,X_2,\dots,X_n)=0$ for any choice of $X_i$.
\end{Proposition}

\subsection{Boolean cumulants}
\label{sec:GeneralCumulantExamples:Boolean}
\emph{Boolean convolution} of measures was studied in 
\cite{SpeicherWoroudi:1997:boolean,BozejkoSpeicher:1991:psi}.
It comes from the so-called \emph{regular free product} of states on
free products of groups \cite{Bozejko:1987:uniformly}.
\begin{Definition}
  Let $(\alg{A}_i=\IC[X_i], \phi_i)$ be polynomial algebras.
  The regular free product of the states $\phi=\bigstar\phi_i$ is
  the state on the unital free product of the algebras
  $\bigstar\alg{A}_i$ which is given by the rule
  $$
  \phi(X_{i_1}^{k_1}
       X_{i_2}^{k_2}
       \dotsm
       X_{i_n}^{k_n})
  = \phi_{i_1}(X_{i_1}^{k_1})
    \,
    \phi_{i_2}(X_{i_2}^{k_2})
    \dotsm
    \phi_{i_n}(X_{i_n}^{k_n})
  $$
  if $i_j\ne i_{j+1}$ and $k_j> 0$.
\end{Definition}
This is a special case of conditional free products considered in
section~\ref{subsec:ConditionalFreeCumulants} below.
The partitions of relevance here are the interval partitions
considered first by von Waldenfels in
\cite{vonWaldenfels:1973:approach,vonWaldenfels:1975:interval}.
\begin{Proposition}
  If $\pi\in\Pi_n$ is an interval partition,
  then the partitioned expectation factorizes:
  $$
  \phi_\pi(X_1 X_2\cdots X_n)
  = \prod_j \phi(\arrowprod_{i\in\pi_j} X_i)
  $$
\end{Proposition}
More generally, the partitioned expectations factorize along the
irreducible components of the partition.
The cumulants also factorize for interval partitions and moreover they vanish
for any other partition.
\begin{Proposition}
  If $\pi\in\Pi_n\setminus I_n$ then
  $K^\exchB_\pi(X_1, X_2, \dots, X_n)$ vanishes for any choice
  of $X_i$.
\end{Proposition}
\begin{proof}
  If there is a block which is not an interval,
  then it is sliced into at least two parts. We have therefore
  an alternating word in which 
  by Lemma~\ref{lem:GeneralCumulants:subwordexpectationvanishes}
  one (even two) of the factors has zero
  expectation and therefore the expectation of the whole word
  vanishes.
\end{proof}
Thus we can write
$$
\phi_\pi(X_1,X_2,\dots,X_n)
= \sum_{\substack{\sigma\leq\pi\\ \sigma\in I_n}}
   K^\exchB_\sigma(X_1,X_2,\dots,X_n)
$$
and we can apply M\"obius inversion on the lattice
of interval partitions and get the formula
$$
K^\exchB_n(X_1,X_2,\dots,X_n)
= \sum_{\pi\in I_n}
   \phi_\pi(X_1,X_2,\dots,X_n)\,\mu_I(\pi,\hat1_n)
.
$$
The name \emph{boolean cumulants} stems from the fact
that the the set of interval partitions of order $n$ forms a lattice
$I_n$ which is isomorphic to the (boolean) lattice of subsets of the set
$\{1,2,\dots,n-1\}$. There is an obvious antiisomorphism which takes
an interval partition to the set of the endpoints of its blocks,
except the last one which is redundant. Then take the antiisomorphism
of the boolean lattice which consists of taking complements.

\subsection{Conditional free cumulants}
\label{subsec:ConditionalFreeCumulants}

A free product of algebras with pairs of states was defined in
\cite{BozejkoSpeicher:1991:psi}%
\paper{Bozejko/Speicher}{}{}
and \cite{BozejkoLeinertSpeicher:1996:convolution}%
\paper{Bozejko/Leinert/Speicher}{}{},
generalizing both free and boolean free product.
Let $\alg{A}_i$ be algebras with states $\phi_i$, $\psi_i$
($\phi$ may also be operator-valued, see~\cite{Mlotkowski:1997:operator}).
On the free product $\alg{A} = \bigstar\alg{A}_i$ let
$\psi=\bigstar \psi_i$
be the free product state and define $\phi$ be the condition
$$
\phi(X_1X_2\dotsm X_n) = \phi_{i_1}(X_1)\, \phi_{i_2}(X_2)\dotsm \phi_{i_n}(X_n)
$$
whenever $X_j\in\alg{A}_{i_j}$, $i_j\ne i_{j+1}$ and $\psi_{i_j}(X_j)=0$.
The resulting noncommutative probability space is called the
conditional free product of $(\alg{A}_i,\phi_i,\psi_i)$ and denoted
\begin{equation}
  \label{eq:GeneralCumulantExamples:CondFreeProduct}
  (\alg{A},\phi,\psi) = \bigstar(\alg{A}_i,\phi_i,\psi_i)
  .
\end{equation}
\begin{Theorem}[{\cite{BozejkoLeinertSpeicher:1996:convolution}\paper[Thm.~2.2]{BLS}{}{}}]
  If all $\phi_i$, $\psi_i$ are states, then $\phi$ is a state.
\end{Theorem}

One can show that the conditional free product is associative,
that is
$$
( (\alg{A}_1,\phi_1,\psi_1) *  (\alg{A}_2,\phi_2,\psi_2) )* (\alg{A}_3,\phi_3,\psi_3)
=
 (\alg{A}_1,\phi_1,\psi_1) * ( (\alg{A}_2,\phi_2,\psi_2) * (\alg{A}_3,\phi_3,\psi_3) )
.
$$
It follows that the conditional free product
$$
\bigstar_{i\in\IN}(\alg{A},\phi,\psi)
$$
of infinitely many 
copies~$(\alg{A}^{(i)},\phi^{(i)},\psi^{(i)})$ 
of~\eqref{eq:GeneralCumulantExamples:CondFreeProduct}
gives rise to an exchangeability system $\exchCF$ such that
the free factors $\alg{A}^{(i)}$ are $\exchCF$-exchangeable
and moreover the subalgebras $\alg{A}_i\subseteq \alg{A}$ are independent
in the sense of definition~\ref{def:GeneralCumulants:independence}.
We can therefore proceed to compute cumulants.
It turns out that crossing cumulants vanish, just as in the free case:

We can decompose $\alg{A}_i=\IC I\oplus\bub{\alg{A}_i}$ where
$\bub{\alg{A}_i}=\ker \psi_i$ and
$$
\alg{A} = \IC I
          \oplus
          \bigoplus_{i_1\ne i_2\ne\dots}
           \bub{\alg{A}_{i_1}}
           \bub{\alg{A}_{i_2}}
           \dotsm
           \bub{\alg{A}_{i_n}}
$$
\begin{Lemma}[{\cite{BozejkoLeinertSpeicher:1996:convolution}\paper[Lemma 2.1]{BLS}{}{}}]
  Let 
  $X_1= S_1 S_2\dotsm S_n
      \in \bub{\alg{A}_{i_1}} 
          \bub{\alg{A}_{i_2}}
          \dotsm
          \bub{\alg{A}_{i_m}}$ with $i_k\ne i_{k+1}$,
  and let $X_2=T_1 T_2\dotsm T_n
          \in \bub{\alg{A}_{j_1}}
          \bub{\alg{A}_{j_2}}
          \dotsm
          \bub{\alg{A}_{j_n}}$ with $j_k\ne j_{k+1}$.
  \begin{enumerate}
   \item if $i_1\ne j_1$ then
    \begin{equation}
      \label{eq:GeneralCumulantExamples:ConditionalStochasticIndependence}
      \phi(X_1^*X_2) = \phi(X_1^*)\,\phi(X_2)
    \end{equation}
   \item if $Y\in\alg{A}_i$ with $i\ne i_1,j_1$ then
    $$
    \phi(X_1^* Y X_2)
    =\psi_i(Y)\,(\phi(X_1^*X_2) - \phi(X_1^*)\,\phi(X_2))
     + \phi_i(Y) \,\phi(X_1^*)\,\phi(X_2)
    $$
  \end{enumerate}
\end{Lemma}
We will need the following simpler version only.
\begin{Lemma}
  \label{lem:GeneralCumulants:ConditionalFreeFactorization}
  \begin{enumerate}
   \item If $X$, $Y$ are c-free then $\phi(XY)=\phi(X)\,\phi(Y)$.
   \item If $\{X_1,X_2\}$ and $Y$ are c-free, then
    $$
    \phi(X_1YX_2)
    = \phi(X_1)\,\phi(Y)\,\phi(X_2)
      +
      \psi(Y)\,
      (\phi(X_1X_2) - \phi(X_1)\,\phi(X_2))
    $$
  \end{enumerate}
\end{Lemma}
\begin{proof}
  Denote $\xi=\psi(X)$, $\bar{X}=\phi(X)$, $\bub{X} = X-\psi(X)=X-\xi$, etc.
  Then we have
  \begin{align*}
    \phi(XY)
    &= \phi((\bub{X}+\xi)(\bub{Y}+\eta)) \\
    &= \phi(\bub{X})\,\phi(\bub{Y})
       +
       \xi\,\phi(\bub{\eta}) + \phi(\bub{X})\,\eta + \xi\eta\\
    &= (\bar{X}-\xi)(\bar{Y}-\eta)
       +
       \xi\,(\bar{Y}-\eta)
       +
       (\bar{X}-\xi)\,\eta
       +
       \xi\eta\\
    &= \bar{X}\bar{Y}
  \end{align*}
  For the second part,
  \begin{align*}
    \phi(X_1YX_2)
    &= \phi(\bub{X_1} Y X_2) + \xi_1 \bar{Y} \bar{X_2}\\
    &= \phi(\bub{X_1} Y \bub{X_2})
       +
       \xi_2\, \phi(\bub{X_1})\, \bar{Y}
       +
       \xi_1 \bar{X}_2 \bar{Y} \\
    &= \phi(\bub{X_1})\,\phi(\bub{Y})\,\phi(\bub{X_2})
       +
       \phi(\bub{X_1}\bub{X_2})\,\eta
       +
       \xi_2\,\phi(\bub{X_1})\,\bar{Y}
       +
       \xi_1 \bar{X}_2 \bar{Y}\\
    &= (\bar{X}_1-\xi_1) (\bar{Y}-\eta) (\bar{X}_2 - \xi_2)
       +
       (\phi(X_1X_2)-\bar{X}_1 \xi_2 - \xi_1\bar{X}_2 + \xi_1\xi_2)\eta
\\ &\phantom{=+}\hfill
       +
       (\bar{X}_1-\xi_1) \xi_2 \bar{Y}
       +
       \xi_1\bar{X}_2\bar{Y}\\
    &= \bar{X}_1\bar{X}_2\bar{Y}
       -
       \bar{X}_1 \bar{X}_2 \eta
       -
       \bar{X}_1 \xi_2 \bar{Y}
       +
       \bar{X}_1 \xi_2 \eta
       -
       \xi_1 \bar{X}_2 \bar{Y}
       +
       \xi_1 \bar{X}_2 \eta
       +
       \xi_1 \xi_2 \bar{Y}
       -
       \xi_1 \xi_2 \eta
\\ &\phantom=
       +
       \phi(X_1 X_2) \,\eta
       -
       \bar{X}_1 \xi_2 \eta
       -
       \xi_1 \bar{X}_2 \eta
       +
       \xi_1 \xi_2 \eta
       + 
       \bar{X}_1 \xi_2 \bar{Y}
       -
       \xi_1 \xi_2 \bar{Y}
       +
       \xi_1 \bar{X}_2 \bar{Y} \\
    &= \bar{X}_1 \bar{X}_2 \bar{Y}
       +
       \eta(\phi(X_1 X_2) - \bar{X}_1 \bar{X}_2)
  \end{align*}
\end{proof}

In particular, pyramidal independence does not hold (unless $\phi_i=\psi_i$, i.e.\ 
free independence).
Let us consider interval partitions first.
We can work in the full algebra with $\psi=\bigstar\psi_i$.
\begin{Lemma}
  The partitioned $\phi$-cumulants are multiplicative on irreducible
  components.
  Let $\pi=\pi_1\cup\pi_2\cup\dots\cup\pi_m$ where $\pi_j$ are the
  irreducible components.
  Then
  \begin{multline*}
    K^\exchCF_\pi(X_1,X_2,\dots,X_n)\\
    = K^\exchCF_{\pi_1}(X_1,X_2,\dots,X_{n_1})
      \,
      K^\exchCF_{\pi_2}(X_{n_1+1},X_{n_1+2},\dots,X_{n_2})
      \dotsm
      K^\exchCF_{\pi_m}(X_{n_{m-1}+1},X_{n_1+1},\dots,X_{n})
  \end{multline*}
\end{Lemma}
\begin{proof}
  This is clear for moments from stochastic independence
  (Proposition~\ref{lem:GeneralCumulants:ConditionalFreeFactorization})
  and associativity of the c-free product.
\end{proof}
Although pyramidal independence does not hold for moments,
it holds for cumulants in a modified way, namely
one has to distinguish inner and outer blocks.
\begin{Proposition}
  Let~$\pi$ be an irreducible partition with outer connected component~$\pi_1$
  and inner components~$\pi_2$.
  Then
  $$
  K^\exchCF_\pi(X_1,X_2,\dots,X_n)
  = K^\exchCF_{\pi_1}(X_1,X_2,\dots,X_n)
  \,
  K^{\exchF,\psi}_{\pi_2}(X_1,X_2,\dots,X_n)
  $$
\end{Proposition}
where $K^{\exchF,\psi}_{\pi_2}$ is the free cumulant with respect to $\psi$.
Therefore it vanishes unless $\pi_2$ is noncrossing.
\begin{proof}
  We may assume that the partition is irreducible, i.e.,
  there is one outer connected component and one or more inner connected components.
  For simplicity let us assume that there is only one inner component,
  of length $m-k$.
  Then we have
  $$
  K^\exchCF_\pi(X_1, X_2, \dots, X_n)
  = \frac{1}{\prod_{B\in\pi}\abs{B}}
    \,
    \phi( X_1^{\pi,\omega} X_2^{\pi,\omega} \cdots X_n^{\pi,\omega})
  $$
  where 
  $\{ X_1^{\pi,\omega},
      X_2^{\pi,\omega},
      \dots,
      X_k^{\pi,\omega},
      X_{m+1}^{\pi,\omega},
      \dots,
      X_n^{\pi,\omega}
  \}$
  and
  $\{ X_{k+1}^{\pi,\omega},
      X_{k+2}^{\pi,\omega},
      \dots,
      X_m^{\pi,\omega}
  \}$
  are independent.
  The pyramidal law of Lemma~\ref{lem:GeneralCumulants:ConditionalFreeFactorization}
  allows the following factorization:
  \begin{multline*}
    K^\exchCF_\pi(X_1, X_2, \dots, X_n)\\
    = 
      \frac{1}{\prod_{B\in\pi}\abs{B}}
      \biggl(
        \phi( X_1^{\pi,\omega} \dotsm X_k^{\pi,\omega})
        \,
        \phi( X_{k+1}^{\pi,\omega} \dotsm X_m^{\pi,\omega})
        \,
        \phi( X_{m+1}^{\pi,\omega} \dotsm X_n^{\pi,\omega})
        \\ \phantom{xxxxxxxxxxxxx}
        + \psi( X_{k+1}^{\pi,\omega} \dotsm X_m^{\pi,\omega})
        \left(
          \phi( X_1^{\pi,\omega} \dotsm X_k^{\pi,\omega}
                X_{m+1}^{\pi,\omega} \dotsm X_n^{\pi,\omega})
\right.\\ \left.
        - \phi( X_1^{\pi,\omega} \dotsm X_k^{\pi,\omega})
          \,
          \phi(X_{m+1}^{\pi,\omega} \dotsm X_n^{\pi,\omega})
      \right)
      \biggr)
  \end{multline*}
  Now both $\phi(X_1^{\pi,\omega}\dotsm X_k^{\pi,\omega})$
  and  $\phi(X_{m+1}^{\pi,\omega}\dotsm X_n^{\pi,\omega})$
  vanish because at least one block of $\pi_1$ is split into two
  and Lemma~\ref{lem:GeneralCumulants:subwordexpectationvanishes} applies.
  Therefore we are left with one term
  $$
  \psi( X_{k+1}^{\pi,\omega} \dotsm X_m^{\pi,\omega})
  \,
  \phi( X_1^{\pi,\omega} \dotsm X_k^{\pi,\omega}
        X_{m+1}^{\pi,\omega} \dotsm X_n^{\pi,\omega})
  ,
  $$
  which is equal to the claimed value.
\end{proof}
\begin{Proposition}
  $K^\exchCF_\pi(X_1,X_2,\dots,X_n)$ vanishes unless $\pi$ is noncrossing.
\end{Proposition}
\begin{proof}
  Using the above formulas,
  we can reduce the proof to consider the connected components separately.
  For the inner components we know by 
  Proposition~\ref{prop:GeneralCumulantExamples:free:crossingcumvanish}
  that crossing free cumulants vanish.
  Thus it is enough to consider a connected partition $\pi$ with at least two blocks.
  In this case no block of $\pi$ is an interval, therefore
  $$
  K^\exchCF_\pi(X_1,X_2,\dots,X_n)
  = \phi(X_1^{\pi,\omega} X_2^{\pi,\omega}\dotsm X_n^{\pi,\omega})
  = \phi(Y_1 Y_2 \dotsm Y_m)
  $$
  is the expectation of an alternating word $Y_1 Y_2\cdots Y_m$
  whose letters satisfy $\phi(Y_j)=\psi(Y_j)=0$ by 
  Lemma~\ref{lem:GeneralCumulants:subwordexpectationvanishes}
  and therefore the expectation vanishes.
\end{proof}

\begin{Remark}
  Besides boolean convolution (see
  section~\ref{sec:GeneralCumulantExamples:Boolean}),
  which corresponds to the state $\psi=\delta_0$ on
  the polynomial algebra $\IC[X]$,
  several other choices of $\psi$ have been studied 
  \cite{Bozejko:2001:deformed,KrystekYoshida:2002:deformed,Yoshida:2002:remarks,Yoshida:2002:weight}.
  Some of these can be reduced to the following ``$\Delta$-convolution'' 
  \cite{Bozejko:2001:deformed,Yoshida:2002:weight}:
  Let $\mu$ be a probability measure on the real line and define a state $\phi$ on $\IC[X]$ by
  $\phi(X^k) = \int t^k d\mu(t)$.
  Let $\omega$ be any probability measure with moments $(\omega_n)_{n\in\IN}$.
  Let $\psi = \omega \boxdot \phi$, that is
  $$
  \psi(X^n) = \omega_n\, \phi(X^n)
  .
  $$
  A certain moment-cumulant formula was found 
  in~\cite{Yoshida:2002:weight}, namely
  $$
  \phi(X^n) = \sum_{\pi\in\NC_n} w(\pi)\, \alpha_\pi^\Delta(X)
  $$
  where $w(\pi)$ is the products of the lengths of all ``arcs'' of $\pi$.
  However, the term $w(\pi) \alpha_\pi^\Delta$ is different
  from the corresponding term $K^\exchCF_\pi(X)$ in the moment-cumulant formula
  of Proposition~\ref{prop:GeneralCumulants:K2M}, which
  corresponds to the conditional free cumulants.
\end{Remark}

\subsection{Fermions and graded indendence}
\label{sec:GeneralCumulantExamples:Fermions}
 Tensor independence, free independence and boolean independence
are the only possible notions of independence in a certain natural
axiomatic scheme \cite{Speicher:1997:universal}.
There are however other notions of independence if an additional structure
is imposed on the noncommutative probability space.
One such example is $\IZ_2$-graded independence
\cite{MingoNica:1997:crossings}.
There is generalization \cite{Goodman:2002:Zngraded} to $\IZ_n$-graded independence, which however does not give rise to
interchangeable algebras and thus does not fit in our framework.
\begin{Definition}
  A $\IZ_2$-graded noncommutative probability space
  $(\alg{A},\gamma,\phi)$ consists of a $\IZ_2$-graded algebra 
  $\alg{A}=\alg{A}_+\oplus\alg{A}_-$, a unital linear functional $\phi$
  and a grading automorphism~$\gamma$ of order~$2$ such that $\phi\circ\gamma=\phi$.
  The elements $X$ of $\alg{A}_\pm$ are called \emph{homogeneous}
  and satisfy $\gamma(X)=\pm X=(-1)^{\partial X}X$ for $X\in\alg{A}_\pm$,
  where the \emph{degree} $\partial X$ is defined as
  $\partial X=0$ if $X\in\alg{A}_+$
  and $\partial X=1$ if $X\in\alg{A}_-$.
  A subalgebra of $\alg{A}$ is called \emph{homogeneous} if it is invariant under $\gamma$.
  Homogeneous subalgebras $\alg{A}_1$, $\alg{A}_2$ of $\alg{A}$ are \emph{graded independent}
  if 
  \begin{enumerate}
   \item they \emph{gradedly commute}, i.e., homogeneous elements $X_1\in \alg{A}_1$
    and $X_2\in \alg{A}_2$ satisfy
    $X_1 X_2 = (-1)^{\partial X_1\partial X_2} X_2X_1$.
   \item $\phi(X_1 X_2) = \phi(X_1)\,\phi(X_2)$ for all $X_i\in\alg{A}_i$.
  \end{enumerate}
  It follows that for odd elements $X$ (i.e., $\partial X=1$) the expectation $\phi(X)=0$.
\end{Definition}
Examples of graded independence include Clifford algebras and the rotation algebra
$\alg{A}_{1/2}$, cf.\  \cite{MingoNica:1997:crossings}.

Here we only recall the graded tensor product, which we will use to construct graded
independent copies of a given algebra.

\begin{Definition}
  Let $(\alg{A}_1, \gamma_1, \phi_1)$ and $(\alg{A}_2, \gamma_2, \phi_2)$ be
  graded non-commutative probability spaces.
  Their graded tensor product $(\alg{A}_1\ox_2\alg{A}_2,\gamma,\phi)$ is defined
  as the usual tensor product $\alg{A}_1\ox\alg{A}_2$ with multiplication
  $$
  (X_1\ox X_2)(X_1'\ox X_2') = (-1)^{\partial X_1\partial X_2} X_1X_1'\ox X_2X_2'
  $$
  for homogeneous elements $X_1,X_1'\in\alg{A}_1$ and $X_2,X_2'\in\alg{A}_2$.
  For arbitrary elements the product is defined by bilinear extension.
  If $\alg{A}_1$ and $\alg{A}_2$ are $*$-algebras,
  then we can make $\alg{A}_1\ox_2\alg{A}_2$ into a star algebra with involution
  $$
  (X_1\ox X_2)^* = (-1)^{\partial X_1\partial X_2} X_1^*\ox X_2^*.
  $$
  The expectation functional is as usual $\phi=\phi_1\ox\phi_2$.
\end{Definition}
It can be shown that the graded tensor product is associative.
Moreover, it gives rise to exchangeable algebras.
\begin{Proposition}
  \label{prop:GeneralCumulantExamples:GradedTensorProductisexchangeable}
  Let $\alg{U}$ be the infinite graded tensor product of copies of 
  the graded noncommutative probability space $\alg{A},\gamma,\phi$ and let
  $\alg{A}_k$ be the $k$-th copy of $\alg{A}$ in $\alg{U}$. As usual we denote
  for $X\in\alg{A}$ its image in $A_k$ by $X^{(k)}$.
  Then $(A_k)_{k\in\IN}$ are interchangeable.
\end{Proposition}
\begin{proof}
  By associativity of the graded tensor product and by multilinearity
  it is enough to show that for homogeneous elements $X_i$, $Y_i\in\alg{A}$
  we have
  $$
  \phi(X_1^{(1)} Y_1^{(2)}X_2^{(1)} Y_2^{(2)}\cdots X_n^{(1)} Y_n^{(2)})
  =  \phi(X_1^{(2)} Y_1^{(1)}X_2^{(2)} Y_2^{(1)}\cdots X_n^{(2)} Y_n^{(1)})
  $$
  We proceed by commuting the tensors:
  \begin{align*}
      \phi(X_1^{(1)} Y_1^{(2)}X_2^{(1)} Y_2^{(2)}\cdots X_n^{(1)} Y_n^{(2)})
      &= \phi((X_1\ox I)(I\ox Y_1) (X_2\ox I)(I\ox Y_2) \cdots (X_n\ox I)(I\ox Y_n)) \\
      &= \phi((X_1\ox Y_1) (X_2\ox Y_2) \cdots (X_n\ox Y_n)) \\
      &= (-1)^{\partial Y_1\partial X_2}
         \phi((X_1X_2\ox Y_1Y_2) (X_3\ox Y_3) \cdots (X_n\ox Y_n)) \\
      &= \cdots \\
      &= (-1)^{\sum_{k=2}^n(\partial Y_1+\cdots+\partial Y_{k-1})\partial X_k}
         \phi(X_1X_2\cdots X_n\ox Y_1Y_2\cdots Y_n)\\
      &= (-1)^{\sum_{i<j}\partial Y_i\partial X_j}
         \phi(X_1X_2\cdots X_n)
         \,
         \phi( Y_1Y_2\cdots Y_n).
  \end{align*}
  On the other hand,
  \begin{align*}
    \phi(X_1^{(2)} Y_1^{(1)}X_2^{(2)} Y_2^{(1)}\cdots X_n^{(2)} Y_n^{(1)})
      &= \phi((I\ox X_1)(Y_1\ox I) (I\ox X_2)(Y_2\ox I) \cdots (I\ox X_n)(Y_n\ox I)) \\
      &= (-1)^{\partial X_1\partial Y_1 + \cdots +\partial X_n\partial Y_n}
         \phi((Y_1\ox X_1) (Y_2\ox X_2) \cdots (Y_n\ox X_n)) \\
      &= (-1)^{\sum_i\partial X_i\partial Y_i + \sum_{i<j}\partial X_i\partial Y_j}
         \phi( Y_1Y_2\cdots Y_n)
         \,
         \phi(X_1X_2\cdots X_n).
  \end{align*}
  Now unless both $\sum \partial Y_i$ and $\sum\partial X_i$ are even,
  the expectations vanish, so we may assume that they are even.
  In this case, the signs are equal, as their product is $1$:
  $$
  (-1)^{\sum_{i>j}\partial X_i\partial Y_j}
  \cdot
  (-1)^{\sum_i\partial X_i\partial Y_i + \sum_{i<j}\partial X_i\partial Y_j}
  = (-1)^{(\sum \partial X_i)(\sum \partial Y_j)}
  = 1
  $$
\end{proof}
\begin{Proposition}
  Graded independence coincides with $\exchGrad$-independence induced,
  via Definition~\ref{def:GeneralCumulants:independence}, by the exchangeability system~$\exchGrad$ of
  Proposition~\ref{prop:GeneralCumulantExamples:GradedTensorProductisexchangeable}.
\end{Proposition}
\begin{Proposition}
  Pyramidal independence holds:
  If $\{X_1, X_2\}$ and $Y$ are graded independent, then
  $$
  \phi(X_1 Y X_2) = \phi(X_1 X_2) \,\phi(Y)
  $$
\end{Proposition}
\begin{proof}
  We may assume that all random variables involved are homogeneous,
  then we haven
  $$
  \phi(X_1YX_2) = (-1)^{\partial X_2\partial Y} \phi(X_1X_2)\,\phi(Y)
  $$
  and $\phi(Y)=0$ unless $\partial Y=0$.
\end{proof}

It follows by Proposition~\ref{prop:GeneralCumulants:PyramidalIndependence}
that moments and cumulants factorize along connected components.
For general partitions the partitioned moments and cumulants also factorize,
but with a weight counting the number of a certain kind of crossings.
\begin{Definition}
  Let $\pi=\{B_1,\dots,B_p\}\in\Pi_n$ be a partition
  with $\min{B_i}<\min{B_j}$ for $i<j$.
  For two blocks $A$ and $B$ of $\pi$ with $\min(A)<\min(B)$ denote
  $$
  c_0(A,B) = \#\{(a,b) : a\in A, b\in B, \min(B)<a<b\}
  $$
  the \emph{(left) reduced number of crossings} of these blocks;
  the total number of reduced crossings of $\pi$ is
  $$
  c_0(\pi) = \sum_{i<j} c_0(B_i,B_j)
  $$
\end{Definition}
\begin{Proposition}
  For homogeneous elements of degree $1$ the partitioned moment and cumulant is
  \begin{align*}
    \phi_\pi(X_1,X_2,\dots,X_n)
    &= (-1)^{c_0(\pi)} \prod_{B\in\pi} \phi_B(\prod_{i\in b} X_i)\\
    K^\exchGrad_\pi(X_1,X_2,\dots,X_n)
    &= (-1)^{c_0(\pi)} \prod_{B\in\pi} K^\exchGrad_B(X_i : i\in b)\\
  \end{align*}
\end{Proposition}
\begin{proof}
  By associativity it is enough to prove that for noncommutative random variables
  $X_1,X_2,\dots,X_n\in\alg{A}$ 
  and a partition $\pi$ consisting of two blocks $B_1$ and $B_2$ 
  with $\min(B_1)=1<\min(B_2)$
  $$
  \phi(X_1^{(\pi(1))} X_2^{(\pi(2))}\cdots X_n^{(\pi(n))})
  = (-1)^{c_0(B_1,B_2)}
    \phi(\prod_{i\in B_1} X_i)
    \,
    \phi(\prod_{i\in B_2} X_i).
  $$
  To do this, we determine the effect of commuting
  $X_m^{(1)}$ to the left of an $X_l^{(2)}$,
  where $k$, $l$ and $m$ are such that
  $$
  X_1^{(\pi(1))} X_2^{(\pi(2))}\cdots X_n^{(\pi(n))}
  = X_1^{(1)} X_2^{(1)}\cdots X_k^{(1)}
    X_{k+1}^{(2)} X_{k+2}^{(2)}\cdots X_l^{(2)}
    X_{l+1}^{(1)} X_{l+2}^{(\pi(l+2))}\cdots X_n^{(\pi(n))}
  $$
  On the one hand,
  \begin{multline*}
  X_1^{(1)} X_2^{(1)}\cdots X_k^{(1)}
  X_{k+1}^{(2)} X_{k+2}^{(2)}\cdots X_l^{(2)}
  X_{l+1}^{(1)} X_{l+2}^{(\pi(l+2))}\cdots X_n^{(\pi(n))}\\
  = - X_1^{(1)} X_2^{(1)}\cdots X_k^{(1)}
      X_{k+1}^{(2)} X_{k+2}^{(2)}\cdots X_{l-1}^{(2)}
      X_{l+1}^{(1)} X_l^{(2)} X_{l+2}^{(\pi(l+2))}\cdots X_n^{(\pi(n))}
  \end{multline*}
  on the other hand, the corresponding new partition $\pi'=\{A',B'\}$ has
  $$
  c_0(A',B') =
  \begin{cases}
    c_0(A,B) + 1            & \text{if $l-k\ge2$}\\
    c_0(A,B) - (\abs{B}-1)    & \text{if $l=k+1$}\\
  \end{cases}
  $$
  crossings. Unless $\abs{B}$ is even, the expectation vanishes anyway,
  therefore in any case we have $(-1)^{c_0(A',B')} = -(-1)^{c_0(A,B)}$.
  Repeating this step until the partition becomes uncrossing finishes the proof.
\end{proof}

\subsection{Noncrossing cumulants of type $B$}

Recently \cite{BianeGoodmanNica:2002:noncrossing}
there has  been introduced a framework for noncrossing cumulants of type $B$
which were defined in \cite{Reiner:1997:noncrossing}.
We are indebted to A.~Nica for explaining the model to us.

Roughly speaking the setup is as follows.
Let $\alg{A}$ be an algebra and $V$ and $\alg{A}$-bimodule,
define a multiplication on $\alg{A}\times V$ induced by the matrix multiplication
$$
\begin{bmatrix}
  a & \xi \\
  0 & a
\end{bmatrix}
\cdot
\begin{bmatrix}
  a' & \xi' \\
  0  & a'
\end{bmatrix}
=
\begin{bmatrix}
  aa' & a\cdot\xi'+\xi\cdot a' \\
  0   & aa'
\end{bmatrix}
$$
i.e., we define the multiplication of a pair
$$
(a,\xi)\cdot(a',\xi') = (aa',a\cdot\xi'+\xi\cdot a')
$$
A similar multiplication is defined on $\IC^2$:
\begin{equation}
  \label{eq:GeneralCumulantExamples:typeB:Cmult}
  (\alpha,\beta)\cdot (\alpha',\beta') = (\alpha\alpha',\alpha\beta'+\alpha'\beta)
\end{equation}
For a given pair of functionals $\phi:\alg{A}\to\IC$ and $f:V\to\IC$
we define the $\IC^2$-valued expectation map
\begin{align*}
  \psi: \alg{A}\times V &\to \IC^2\\
      (a,\xi) &\mapsto (\phi(a),f(\xi))
\end{align*}

\begin{Definition}
  Let $\alg{A}_i$ be subalgebras of $\alg{A}$ and $V_i$ subspaces of $V$
  s.t. $V_i$ is invariant under the action of $\alg{A}_i$.
  The family $(\alg{A}_i,V_i)$ is \emph{free of type B} if
  \begin{enumerate}[(i)]
   \item $\alg{A}_i$ are free
   \item Whenever $a_j\in\alg{A}_{i_{-j}}$, $\xi\in \alg{A}_{i_0}$
     and $b_j\in\alg{A}_{i_j}$ with $i_j\ne i_{j+1}$, $\phi(a_j)=\phi(b_j)=0$
     we have
     \begin{equation}
       \label{eq:GeneralCumulantExamples:typeB:freedef2}
       f(a_m a_{n-1} \cdots a_1 \cdot \xi \cdot b_1 b_2\cdots b_n)
       = \begin{cases}
           0 & m\ne n\\
           \delta_{j_{-1} j_1}
           \delta_{j_{-2} j_2}
           \cdots
           \delta_{j_{-n} j_n}
           \phi(a_1 b_1)
           \phi(a_2 b_2)
           \cdots
           \phi(a_n b_n)
           f(\xi)
         \end{cases}
     \end{equation}
  \end{enumerate}
\end{Definition}
There is also a natural free product construction associated to this notion.

\begin{Lemma}
  Pyramidal independence holds.
  Assume that $\{(a,\xi), (a',\xi')\}$ and $(b,\eta)$ are free of type B.
  Then
  $$
  \psi( (a,\xi)\cdot(b,\eta)\cdot(a',\xi') )
  =
  \psi( (aa',\xi\xi') )\cdot \psi( (b,\eta) )
  $$
\end{Lemma}
\begin{proof}
  By definition we have 
  $$
  \psi( (a,\xi)\cdot(b,\eta)\cdot(a',\xi') )
  = (\phi(aba'), f(\xi\cdot ba' + a\cdot\eta\cdot a' + ab\xi'))
  $$
  Pyramidal independence holds in free probability and therefore the first component
  is clearly $\phi(aba') = \phi(aa')\,\phi(b)$.
  The other terms are
  \begin{align*}
    f(\xi\cdot ba')
    &= f( \xi\cdot (\bub{b} + \phi(b)) \, (\bub{a'} + \phi(a'))\\
    &= f( \xi\cdot \bub{b}\bub{a'} )
       +
       f( \xi\cdot \bub{b}) \, \phi(a')
       + 
       f( \xi\cdot \bub{a'}) \, \phi(b) \, f( \xi) \, \phi(b) \, \phi(a') \\
    &= f(\xi\cdot a')\,\phi(b)\\
    f(a\cdot\eta a')
    &= f( (\bub{a}+\phi(a)) \cdot \eta\cdot  (\bub{a'}+\phi(a')) ) \\
    &= f( \bub{a}\cdot\eta\cdot\bub{a'})
       +
       f( \bub{a}\cdot\eta)\,\phi(a')
       +
       \phi(a')\,f( \eta\cdot\bub{a})
       +
       \phi(a')\,f( \eta) \, \phi(a') \\
    &= \phi(\bub{a}\bub{a'})\,f(\eta) + \phi(a)\,f(\eta)\,\phi(a')\\
    &= \phi(aa')\,f(\eta)\\
    f(ab\cdot\xi')
    &= f(a\xi')\,\phi(b)
  \end{align*}
  On the other hand (noting that $\IC^2$ ist commutative with the multiplication
  \eqref{eq:GeneralCumulantExamples:typeB:Cmult}) we have
  \begin{align*}
    \psi( (a,\xi)\cdot(a',\xi') ) \cdot \psi( (b,\eta) )
    &= ( \phi(aa'), f(a\cdot\xi') + f(\xi\cdot a')) \cdot (\phi(b), f(\eta)) \\
    &= ( \phi(aa')\,\phi(b), \phi(aa')f(\eta) + \phi(b)\,(f(a\cdot\xi') + f(\xi\cdot a')))
  \end{align*}
  and this coincides with the value above.
\end{proof}

\begin{Proposition}
  Crossing cumulants vanish.
\end{Proposition}

\begin{proof}
  By Corollary~\ref{prop:GeneralCumulants:PyramidalIndependence} the cumulants factor along
  connected components. It is therefore enough to compute them for connected partitions.
  Assume that $\pi$ is connected and has at least $2$ blocks.
  Consider
  \begin{multline*}
    \psi( (a_1,\xi_1)^{(\pi,\omega)}
          \cdot
          (a_2,\xi_2)^{(\pi,\omega)}
          \cdots
          (a_n,\xi_n)^{(\pi,\omega)}
        ) \\
    = ( \phi(a_1^{(\pi,\omega)}
              a_2^{(\pi,\omega)}
              \cdots
              a_n^{(\pi,\omega)}
             ),
         f( \xi_1^{(\pi,\omega)} \cdot a_2^{(\pi,\omega)}\cdots a_n^{(\pi,\omega)}
            +
            a_1^{(\pi,\omega)}\cdot\xi_2^{(\pi,\omega)} \cdots a_n^{(\pi,\omega)}
            +
            a_1^{(\pi,\omega)} \cdots a_{n-1}^{(\pi,\omega)}\cdot \xi_n^{(\pi,\omega)}))
  \end{multline*}
  The first component vanishes by
  Proposition~\ref{prop:GeneralCumulantExamples:free:crossingcumvanish};
  for the second component consider the summand
  $$
  f(a_1^{(\pi,\omega)}
    \cdots 
    a_{k-1}^{(\pi,\omega)}
    \cdot
    \xi_k^{(\pi,\omega)}
    \cdot
    a_{k+1}^{(\pi,\omega)}
    \cdots
    a_n^{(\pi,\omega)}
   )
  $$
  assume that $\{j,j+1,\dots,k-1,k,k+1,\dots,l\}$ is the maximal interval
  containing $k$ and contained in a block of $\pi$
  by assumption on $\pi$ this interval is not a block of $\pi$ and therefore
  by Lemma~\ref{lem:GeneralCumulants:subwordexpectationvanishes} we have
  $$
  f(a_j^{(\pi,\omega)}
    \cdots
    a_{k-1}^{(\pi,\omega)}
    \cdot
    \xi_k^{(\pi,\omega)}
    \cdot
    a_{k+1}^{(\pi,\omega)}
    \cdots
    a_l^{(\pi,\omega)}
    )
  = 0
  .
  $$
  Now by \eqref{eq:GeneralCumulantExamples:typeB:freedef2}
  it follows that
  $$
  f(a_1^{(\pi,\omega)}
    \cdots 
    a_{k-1}^{(\pi,\omega)}
    \cdot
    \xi_k^{(\pi,\omega)}
    \cdot
    a_{k+1}^{(\pi,\omega)}
    \cdots
    a_n^{(\pi,\omega)}
   )
   = 0
  $$
\end{proof}

To summarize, only noncrossing cumulants contribute and are multiplicative
(w.r.\ to the multiplication \eqref{eq:GeneralCumulantExamples:typeB:Cmult}).
They are expressed in terms of so-called noncrossing cumulants of type B in
\cite{BianeGoodmanNica:2002:noncrossing}.

\subsection{Further examples}
Other interesting examples can be constructed from relatively free groups.
Take any relatively free group on infinitely many generators, 
for instance the free nilpotent group $G=\langle g_1,g_2,\dots | [x,[y,z]]=1\rangle$
or the free metabelian group $G=\langle g_1,g_2,\dots | [[v,w],[x,y]]=1\rangle$
(or, even more daring, the Burnside group $G=\langle g_1,g_2,\dots | x^n=1 \rangle$).
Then the group algebra of such a group contains an interchangeable family of copies
of the group algebra of $\mathbf{Z}$ or the free nilpotent (metabelian) group
on a fixed finite number $N$ of generators.
It seems to be an interesting problem to determine the relevant partition
lattice for these groups, that is, to sieve out those partitions, for which the
cumulants always vanish. We  will deal with this ``free nilpotent
probability'' and ``free metabelian probability'' elsewhere.


\bibliography{GeneralCumulants}

\providecommand{\bysame}{\leavevmode\hbox to3em{\hrulefill}\thinspace}
\begin{thebibliography}{VDN92}

\bibitem[Avi82]{Avitzour:1982:free}
Avitzour, D., \emph{Free products of ${C}\sp{\ast} $-algebras}, Trans. Amer.
  Math. Soc. \textbf{271} (1982), 423--435.

\bibitem[BGN03]{BianeGoodmanNica:2002:noncrossing}
Biane, P., Goodman, F., and Nica, A., \emph{Non-crossing cumulants of type
  {B}}, Trans. Amer. Math. Soc. \textbf{355} (2003), 2263--2303,
  arXiv:math.OA/0206167.

\bibitem[BLS96]{BozejkoLeinertSpeicher:1996:convolution}
Bo{\.z}ejko, M., Leinert, M., and Speicher, R., \emph{Convolution and limit
  theorems for conditionally free random variables}, Pacific J. Math.
  \textbf{175} (1996), 357--388.

\bibitem[Bo{\.z}87]{Bozejko:1987:uniformly}
Bo{\.z}ejko, M., \emph{Uniformly bounded representations of free groups}, J.
  Reine Angew. Math. \textbf{377} (1987), 170--186.

\bibitem[Bo{\.z}01]{Bozejko:2001:deformed}
Bo{\.z}ejko, M., \emph{Deformed free probability of {V}oiculescu},
  S\=urikaisekikenky\=usho K\=oky\=uroku (2001), no.~1227, 96--113, Infinite
  dimensional analysis and quantum probability theory (Japanese) (Kyoto, 2000).

\bibitem[Bri69]{Brillinger:1969:calculation}
Brillinger, D., \emph{The calculation of cumulants via conditioning}, Ann.
  Inst. Statist. Math. \textbf{21} (1969), 375--390.

\bibitem[BS91]{BozejkoSpeicher:1991:psi}
Bo{\.z}ejko, M., and Speicher, R., \emph{$\psi$-independent and symmetrized
  white noises}, Quantum probability \& related topics, World Sci. Publishing,
  River Edge, NJ, 1991, pp.~219--236.

\bibitem[BS96]{BozejkoSpeicher:1996:interpolations}
Bo{\.z}ejko, M., and Speicher, R., \emph{Interpolations between bosonic and
  fermionic relations given by generalized {B}rownian motions}, Math. Z.
  \textbf{222} (1996), 135--159.

\bibitem[CD99]{Cabanal-Duvillard:1999:noncrossing}
Cabanal-Duvillard, T., \emph{Non-crossing partitions for conditional freeness},
  Preprint, 1999.

\bibitem[CT78]{ChowTeicher:1978:probability}
Chow, Y.~S., and Teicher, H., \emph{Probability theory}, Springer-Verlag, New
  York, 1978.

\bibitem[DRS72]{DoubiletRotaStanley:1972:foundationsVI}
Doubilet, P., Rota, G.-C., and Stanley, R., \emph{On the foundations of
  combinatorial theory. {V}{I}. {T}he idea of generating function}, Proceedings
  of the Sixth Berkeley Symposium on Mathematical Statistics and Probability,
  Vol. II: Probability theory, Univ. California Press, Berkeley, Calif., 1972,
  pp.~267--318.

\bibitem[Fis29]{Fisher:1929:moments}
Fisher, R., \emph{Moments and product moments of sampling distributions}, Proc.
  Lond. Math. Soc. Series 2 \textbf{30} (1929), 199--238, Reprinted as paper 74
  in \emph{Collected Papers of R.A. Fisher}, vol. 2, (ed. J.H. Bennett. Univ.
  of Adelaide Press, 1972, 351--354.

\bibitem[Goo75]{Good:1975:new}
Good, I.~J., \emph{A new formula for cumulants}, Math. Proc. Cambridge Philos.
  Soc. \textbf{78} (1975), 333--337.

\bibitem[Goo77]{Good:1977:new}
Good, I.~J., \emph{A new formula for $k$-statistics}, The Ann. of Statist.
  \textbf{5} (1977), 224--228.

\bibitem[Goo02]{Goodman:2002:Zngraded}
Goodman, F.~M., \emph{{$\mathbb Z_n$--graded Independence}}, 2002,
  arXiv:math.OA/0206296.

\bibitem[GS01]{BenGhorbalSchurmann:2001:algebraic}
Ghorbal, A.~B., and Sch{\"u}rmann, M., \emph{On the algebraic foundations of
  non-commutative probability theory}, Preprint, 2001.

\bibitem[Hal00]{Hald:2000:early}
Hald, A., \emph{The early history of the cumulants and the {G}ram--{C}harlier
  series}, Internat. Statist. Rev. \textbf{68} (2000), 137--153.

\bibitem[Heg85]{Hegerfeldt:1985:noncommutative}
Hegerfeldt, G.~C., \emph{Noncommutative analogs of probabilistic notions and
  results}, J. Funct. Anal. \textbf{64} (1985), 436--456.

\bibitem[HP00]{HiaiPetz:2000:semicircle}
Hiai, F., and Petz, D., \emph{The semicircle law, free random variables and
  entropy}, American Mathematical Society, Providence, RI, 2000.

\bibitem[Kin78]{Kingman:1978:uses}
Kingman, J. F.~C., \emph{Uses of exchangeability}, Ann. Probability \textbf{6}
  (1978), 183--197.

\bibitem[Kre72]{Kreweras:1972:partitions}
Kreweras, G., \emph{Sur les partitions non crois\'ees d'un cycle}, Discrete
  Math. \textbf{1} (1972), 333--350.

\bibitem[KS00]{KrawczykSpeicher:2000:combinatorics}
Krawczyk, B., and Speicher, R., \emph{Combinatorics of free cumulants}, J.
  Combin. Theory Ser. A \textbf{90} (2000), 267--292.

\bibitem[KY02]{KrystekYoshida:2002:deformed}
Krystek, A., and Yoshida, H., \emph{Deformed {N}arayana number arising from the
  $r$-free convolution}, Preprint, 2002.

\bibitem[Leh02]{Lehner:2002:connected}
Lehner, F., \emph{Free cumulants and enumeration of connected partitions},
  European J. Combin. \textbf{23} (2002), 1025--1031.

\bibitem[Leh03a]{Lehner:2002:Cumulants2}
Lehner, F., \emph{Cumulants in noncommutative probability theory {II}.
  {G}eneralized {G}aussian random variables}, Probab. Theory Related Fields
  \textbf{127} (2003), 407--422, arXiv:math.CO/0210443.

\bibitem[Leh03b]{Lehner:2002:Cumulants3}
Lehner, F., \emph{Cumulants in noncommutative probability theory {III}.
  {C}reation- and annihilation operators on {F}ock spaces}, Preprint,
  arXiv:math.CO/0210444, 2003.

\bibitem[Leh04]{Lehner:2004:cumulants4}
Lehner, F., \emph{Cumulants in noncommutative probability theory {IV}.
  {N}oncrossing cumulants: {D}e {F}inetti's theorem, {$L^p$}-inequalities and
  {B}rillinger's formula}, Preprint, 2004.

\bibitem[LS59]{LeonovShiryaev:1959:method}
Leonov, V.~P., and Shiryaev, A.~N., \emph{On a method of calculation of
  semi-invariants}, Theor. Prob. Appl. \textbf{4} (1959), 319--328.

\bibitem[Mat99]{Mattner:1999:what}
Mattner, L., \emph{What are cumulants?}, Doc. Math. \textbf{4} (1999),
  601--622.

\bibitem[M{\l}o02]{Mlotkowski:1997:operator}
M{\l}otkowski, W., \emph{Operator-valued version of conditionally free
  product}, Studia Math. \textbf{153} (2002), 13--30.

\bibitem[MN97]{MingoNica:1997:crossings}
Mingo, J.~A., and Nica, A., \emph{Crossings of set-partitions and addition of
  graded-independent random variables}, Internat. J. Math. \textbf{8} (1997),
  645--664.

\bibitem[NS00]{NicaSpeicher:2000:combinatorics}
Nica, A., and Speicher, R., \emph{Combinatorics of free probability theory},
  Lecture notes from a course at the IHP, Paris 1999, 2000.

\bibitem[NSS02]{NicaShlyakhtenkoSpeicher:2001:Rcyclic}
Nica, A., Shlyakhtenko, D., and Speicher, R., \emph{${R}$-cyclic families of
  matrices in free probability}, J. Funct. Anal. \textbf{188} (2002), 227--271.

\bibitem[Rei97]{Reiner:1997:noncrossing}
Reiner, V., \emph{Non-crossing partitions for classical reflection groups},
  Discrete Math. \textbf{177} (1997), 195--222.

\bibitem[Rot64]{Rota:1964:foundationsI}
Rota, G.-C., \emph{On the foundations of combinatorial theory. {I}. {T}heory of
  {M}{\"o}bius functions}, Z. Wahrscheinlichkeitstheorie und Verw. Gebiete
  \textbf{2} (1964), 340--368.

\bibitem[Sch47]{Schutzenberger:1947:certains}
Sch\"utzenberger, M.-P., \emph{Sur certains param\`etres caract\'eristiques des
  syst\`emes d'\'ev\'enements compatibles et d\'ependants et leur application
  au calcul des cumulants de la r\'ep\'etition}, C. R. Acad. Sci. Paris
  \textbf{225} (1947), 277--278.

\bibitem[Sch95]{Schurmann:1995:direct}
Sch{\"u}rmann, M., \emph{Direct sums of tensor products and non-commutative
  independence}, J. Funct. Anal. \textbf{133} (1995), 1--9.

\bibitem[Spe83]{Speed:1983:cumulantsI}
Speed, T.~P., \emph{Cumulants and partition lattices}, Austral. J. Statist.
  \textbf{25} (1983), 378--388.

\bibitem[Spe94]{Speicher:1994:multiplicative}
Speicher, R., \emph{Multiplicative functions on the lattice of noncrossing
  partitions and free convolution}, Math. Ann. \textbf{298} (1994), 611--628.

\bibitem[Spe97]{Speicher:1997:universal}
Speicher, R., \emph{On universal products}, Free probability theory (Waterloo,
  ON, 1995), Amer. Math. Soc., Providence, RI, 1997, pp.~257--266.

\bibitem[Spe98]{Speicher:1998:combinatorial}
Speicher, R., \emph{Combinatorial theory of the free product with amalgamation
  and operator-valued free probability theory}, Mem. Amer. Math. Soc.
  \textbf{132} (1998), no.~627, x+88.

\bibitem[Spe00]{Speicher:2000:conceptual}
Speicher, R., \emph{A conceptual proof of a basic result in the combinatorial
  approach to freeness}, Infin. Dimens. Anal. Quantum Probab. Relat. Top.
  \textbf{3} (2000), 213--222.

\bibitem[SW97]{SpeicherWoroudi:1997:boolean}
Speicher, R., and Woroudi, R., \emph{Boolean convolution}, Free probability
  theory (Waterloo, ON, 1995), Amer. Math. Soc., Providence, RI, 1997,
  pp.~267--279.

\bibitem[VDN92]{VDN:1992:free}
Voiculescu, D.~V., Dykema, K.~J., and Nica, A., \emph{Free random variables},
  CRM Lecture Notes Series, vol.~1, American Mathematical Society, Providence,
  RI, 1992.

\bibitem[vLM96]{vanLeeuwenMaassen:1996:obstruction}
van Leeuwen, H., and Maassen, H., \emph{An obstruction for $q$-deformation of
  the convolution product}, J. Phys. A \textbf{29} (1996), 4741--4748.

\bibitem[Voi85]{Voiculescu:1985:symmetries}
Voiculescu, D., \emph{Symmetries of some reduced free product ${C}\sp
  \ast$-algebras}, Operator algebras and their connections with topology and
  ergodic theory (Bu\c steni, 1983), Springer, Berlin, 1985, pp.~556--588.

\bibitem[Voi95]{Voiculescu:1995:operations}
Voiculescu, D., \emph{Operations on certain non-commutative operator-valued
  random variables}, Ast\'erisque (1995), no.~232, 243--275, Recent advances in
  operator algebras (Orl\'eans, 1992).

\bibitem[Voi00]{Voiculescu:2000:lectures}
Voiculescu, D., \emph{Lectures on free probability theory}, Lectures on
  probability theory and statistics (Saint-Flour, 1998), Springer, Berlin,
  2000, pp.~279--349.

\bibitem[vW73]{vonWaldenfels:1973:approach}
von Waldenfels, W., \emph{An approach to the theory of pressure broadening of
  spectral lines}, Probability and information theory, II, Springer, Berlin,
  1973, pp.~19--69. Lecture Notes in Math., Vol. 296.

\bibitem[vW75]{vonWaldenfels:1975:interval}
von Waldenfels, W., \emph{Interval partitions and pair interactions},
  S\'eminaire de Probabilit\'es, IX, Springer, Berlin, 1975, pp.~565--588.
  Lecture Notes in Math., Vol. 465.

\bibitem[Yos02a]{Yoshida:2002:remarks}
Yoshida, H., \emph{Remarks on the $s$-free convolution}, Preprint, 2002.

\bibitem[Yos02b]{Yoshida:2002:weight}
Yoshida, H., \emph{The weight function on non-crossing partitions for the
  ${\Delta}$-convolution}, Preprint, 2002.

\end{thebibliography}
\bibliographystyle{mamsalpha}
\end{document}